\documentclass[stsy,nonblindrev]{informs-stsy}

\OneAndAHalfSpacedXI
\usepackage{natbib}
 \bibpunct[, ]{(}{)}{,}{a}{}{,}%
 \def\bibfont{\small}%
\TheoremsNumberedThrough     % Preferred (Theorem 1, Lemma 1, Theorem 2)
\ECRepeatTheorems

\EquationsNumberedThrough    % Default: (1), (2), ...
\MANUSCRIPTNO{}

\usepackage[letterpaper]{geometry}
\usepackage{authblk}
\usepackage{amsmath,amsfonts}
\usepackage{color,graphicx}
\usepackage[colorlinks,linkcolor=blue,citecolor=cyan,anchorcolor=blue]{hyperref}
\usepackage{enumerate}
\usepackage{wrapfig}
\usepackage{accents}
\usepackage{tikz}
\usetikzlibrary{arrows,
                petri,
                topaths}
\usepackage{multirow}
\usepackage{verbatim}
\usetikzlibrary{shapes,calc,arrows.meta,positioning}
\usetikzlibrary{plotmarks}
\usetikzlibrary{shapes}

\newcommand{\supp}{\text{supp}}
\newcommand{\R}{\mathbb{R}}

\newcommand{\Z}{\mathbb{Z}}
\newcommand{\N}{\mathbb{N}}
\newcommand{\E}{\mathbb{E}}

\newcommand{\lipone}{\text{\rm Lip(1)}}
\newcommand{\dlipone}{\text{\rm dLip(1)}}
\newcommand{\Prob}{\mathbb{P}}

\newcommand{\startproof}{\proof}
\newcommand{\finishproof}{\hfill $\square$\endproof}
\providecommand{\abs}[1]{\left\lvert#1\right\rvert}
\providecommand{\norm}[1]{\lVert#1\rVert}
\usepackage[normalem]{ulem}

\begin{document}
%%%%%%%%%%%%%%%%

% Outcomment only when entries are known. Otherwise leave as is and
%   default values will be used.
%\setcounter{page}{1}
%\VOLUME{00}%
%\NO{0}%
%\MONTH{Xxxxx}% (month or a similar seasonal id)
%\YEAR{0000}% e.g., 2005
%\FIRSTPAGE{000}%
%\LASTPAGE{000}%
%\SHORTYEAR{00}% shortened year (two-digit)
%\ISSUE{0000} %
%\LONGFIRSTPAGE{0001} %
%\DOI{10.1287/xxxx.0000.0000}%

% Author's names for the running heads
% Sample depending on the number of authors;
% \RUNAUTHOR{Jones}
% \RUNAUTHOR{Jones and Wilson}
% \RUNAUTHOR{Jones, Miller, and Wilson}
% \RUNAUTHOR{Jones et al.} % for four or more authors
% Enter authors following the given pattern:
\RUNAUTHOR{Braverman}

% Title or shortened title suitable for running heads. Sample:
% \RUNTITLE{Bundling Information Goods of Decreasing Value}
% Enter the (shortened) title:
\RUNTITLE{Prelimit generator comparison approach}

% Full title. Sample:
% \TITLE{Bundling Information Goods of Decreasing Value}
% Enter the full title:
\TITLE{The prelimit generator comparison approach of Stein's method}

% Block of authors and their affiliations starts here:
% NOTE: Authors with same affiliation, if the order of authors allows,
%   should be entered in ONE field, separated by a comma.
%   \EMAIL field can be repeated if more than one author
\ARTICLEAUTHORS{%
\AUTHOR{Anton Braverman}
\AFF{Kellogg School of Management, Northwestern University, Evanston, IL 60208, \EMAIL{anton.braverman@kellogg.northwestern.edu}} %, \URL{}}
% Enter all authors
} % end of the block

\ABSTRACT{%
This paper uses the generator comparison approach of Stein's method to analyze the gap between  steady-state distributions of Markov chains and diffusion processes. The ``standard'' generator comparison approach starts with the Poisson equation for the diffusion, and the main technical difficulty is to obtain bounds on the derivatives of the solution to the Poisson equation, also known as Stein factor bounds. In this paper we propose starting with the Poisson equation of the Markov chain; we term this the \emph{prelimit approach}. Although one still needs Stein factor bounds, they now correspond to  finite differences of the Markov chain Poisson equation solution rather than the derivatives of the solution to the diffusion Poisson equation. In certain cases, the former are easier to obtain. We use the $M/M/1$ model as a simple working example to illustrate our approach.
% Enter your abstract
}%

% Sample
%\KEYWORDS{deterministic inventory theory; infinite linear programming duality;
%  existence of optimal policies; semi-Markov decision process; cyclic schedule}

% Fill in data. If unknown, outcomment the field

\maketitle
%%%%%%%%%%%%%%%%%%%%%%%%%%%%%%%%%%%%%%%%%%%%%%%%%%%%%%%%%%%%%%%%%%%%%%

% Samples of sectioning (and labeling) in MNSC
% NOTE: (1) \section and \subsection do NOT end with a period
%       (2) \subsubsection and lower need end punctuation
%       (3) capitalization is as shown (title style).
%
%\section{Introduction.}\label{intro} %%1.
%\subsection{Duality and the Classical EOQ Problem.}\label{class-EOQ} %% 1.1.
%\subsection{Outline.}\label{outline1} %% 1.2.
%\subsubsection{Cyclic Schedules for the General Deterministic SMDP.}
%  \label{cyclic-schedules} %% 1.2.1
%\section{Problem Description.}\label{problemdescription} %% 2.

% Text of your paper here

\section{Introduction}
\label{sec:introduction}
Recent years have seen growing use of the generator comparison approach of Stein's method to establish rates of convergence for steady-state diffusion approximations of Markov chains.  One very active area has been the study of queueing and service systems, e.g.,\  \cite{Stol2015,Gurv2014, BravDai2017,BravDaiFeng2016,Ying2016, Ying2017, DaiShi2017, GurvHuan2018, FengShi2018, LiuYing2019, bravgurvhuan2020, Brav2020, BravDaiFang2020}. In the typical setup, one considers a parametric family of continuous-time Markov chains (CTMCs) $\{X(t)\}$ taking values in some discrete state space. This family is often termed the \emph{prelimit} sequence. As the parameters tend to some asymptotic limit, the prelimit sequence converges to a limiting diffusion process $\{Y(t)\}$.  In queueing, for example, the CTMC parameters are usually the arrival rate, number of servers, and service rate, and one common asymptotic regime  is  where the system utilization approaches one, also known as the heavy-traffic regime. To allow for general CTMC families, we assume the CTMC takes values in $\delta \Z^{d} = \{\delta k :\ k \in \Z^{d}\}$ where $\delta > 0$ is a parameter of the CTMC and  the asymptotic regime of interest has $\delta$ converging to zero.  To simplify notation, we omit the dependence of the CTMC on $\delta$ (or any other parameters).  Let  $X$ and $Y$ denote vectors having the stationary distribution of the CTMC and diffusion, respectively. We emphasize that these refer to the stationary distributions and not the stochastic processes $\{X(t)\}$ and $\{Y(t)\}$. The generator approach of Stein's method has been used to study the rates of convergence of $X$ to $Y$ as $\delta \to 0$. The generator approach is attributed to  \cite{Barb1988,Barb1990} and \cite{Gotz1991}, which were the first papers to connect Stein's method to generators of diffusions and CTMCs.

The limiting factor in the generator comparison approach is the curse of dimensionality, because the distance between $X$ and $Y$ depends on the derivatives of the solution to the Poisson equation of the diffusion. In the literature on Stein's method, the Poisson equation is also referred to as the Stein equation. When the diffusion is multidimensional, the Poisson equation is a second-order partial differential equation (PDE), and obtaining derivative bounds, also known as Stein factor bounds, becomes a challenge. The present paper is concerned with expanding the technical toolbox for getting multidimensional Stein factor bounds.  Before discussing our contribution, let us examine this problem in detail.

Recall that $X \in \delta \Z^{d}$ and assume $Y \in \R^{d}$. Let $G_X$ and $G_Y $ be  the infinitesimal generators of the CTMC and diffusion, respectively. Suppose $G_X$ has the form
\begin{align}
G_X f(\delta k) =&\ \sum_{k' \in   \Z^{d} } q_{k,k'} (f(\delta k')-f(\delta k)), \quad k \in  \Z^{d}, \label{eq:generalctmc}
\end{align}
where $q_{k,k'}$ are the transition rates from $\delta k$ to $\delta k'$. Further suppose the diffusion generator has the form
\begin{align*}
G_Y f(x) = \sum_{i=1}^{d} b_i(x) \frac{\partial}{\partial x_i} f(x) + \frac{1}{2} \sum_{i,j=1}^{d} a_{ij}(x) \frac{\partial^2}{\partial x_i \partial x_j} f(x), \quad x \in \R^{d},
\end{align*}
where $f: \R^{d}\to \R$ is a twice continuously differentiable function and $b(x) = (b_1(x), \ldots, b_d(x))$ and $a(x) = (a_{ij}(x))_{i,j=1}^{d}$ are known as the drift and diffusion coefficient, respectively.   

The generator approach works as follows. First, we choose a test function $h^{*}: \R^{d} \to \R$ and consider the Poisson equation  
\begin{align}
G_Y f_{h^{*}}(x) =&\ \E h^{*}(Y) - h^{*}(x), \quad x \in \R^{d}. \label{eq:poisson_diffusion}
\end{align}
We use the star superscript above to emphasize that the functions are defined on all of $\R^{d}$. Given an arbitrary random element $W \in \R^{d}$, one can compare $\E h^{*}(Y)$ to $\E h^{*}(W)$ by taking expected values with respect to $W$ above and attempting to bound the left-hand side. Choosing $W = X$ allows us to leverage the fact that  $\E G_X f_{h^{*}}(X) = 0$ (under some mild conditions), and so 
\begin{align}
\E h^{*}(Y) - \E h^{*}(X) = \E \big( G_Y f_{h^{*}}(X)  - G_X f_{h^{*}}(X) \big). \label{eq:main_difference}
\end{align}
In practice, the chosen $h^{*}(x)$ frequently belongs to 
\begin{align*}
\lipone = \{h^{*}:  \R^{d} \to \R : \abs{h(x)-h(y)} \leq \abs{x-y}, \text{ for all } x,y \in \R^{d} \}.
\end{align*}
This choice is made because functions in $\lipone$ are simple to work with, and because the Wasserstein distance 
\begin{align*}
d_{\lipone}(X,Y) = \sup_{h^{*} \in \lipone} \big| \E h^*(X) - \E h^*(Y) \big|
\end{align*} 
is convergence determining; i.e.,\  convergence in the Wasserstein distance implies convergence in distribution (see, for instance, \cite{GibbSu2002}).  

Bounding the error on the right hand side of \eqref{eq:main_difference}   requires bounds on the derivatives of $f_{h^{*}}(x)$. We refer to these as ``derivative bounds.'' Depending on the transition structure of the CTMC, one may also need to bound certain moments of $X$.  Usually, the approximation $Y$ is such that \eqref{eq:main_difference} converges to zero at a rate of $\delta$, and to prove this it  suffices to bound the second and third derivatives of $f_{h^{*}}(x)$. However, when one seeks approximations $Y$ with convergence rates  faster than $\delta$, as was done in \cite{BravDaiFang2020}, for example, one needs to bound  fourth- and higher-order derivatives.
%In \cite{BravDaiFeng2016}, these are referred to as gradient bounds, and moment bounds, respectively. Another commonly used term for the gradient bounds is Stein factors.

When $d = 1$, the explicit form of $f_{h^{*}}(x)$ is known and can be used to get the derivative bounds via a brute-force approach. When $d > 1$, the Poisson equation is a second-order PDE, and the same kind of brute-force analysis cannot be carried out. Instead, one has to rely on the fact that provided it is finite, 
\begin{align}
f_{h^{*}}(x) = \int_{0}^{\infty} \big( \E_{Y(0)=x} h^{*}(Y(t)) - \E h^{*}(Y) \big) dt, \quad x \in \R^{d}, \label{eq:diffusionpoissonsolution}
\end{align} 
solves the Poisson equation; see any one of  \cite{Barb1990,Gotz1991, Gurv2014,GorhMack2016}  for a proof. We then have that
\begin{align}
&\frac{\partial}{\partial x_i} f_{h^{*}}(x) \approx \frac{f_{h^{*}}(x + \varepsilon e^{(i)}) - f_{h^{*}}(x)}{\varepsilon}  =  \frac{1}{\varepsilon}  \int_{0}^{\infty} \big( \E_{Y(0) = x + \varepsilon e^{(i)}} h^{*}(Y(t)) - \E_{Y(0) = x} h^{*}(Y(t)) \big) dt. \label{eq:first_deriv}
\end{align}
Higher-order derivatives can be accessed similarly. 
There are a few ways to bound  \eqref{eq:first_deriv}. In a handful of cases, the  distribution of $Y(t)$ is known as a function of $Y(0)$, as in \cite{Barb1990,Gotz1991,GanRollRoss2017, GanRoss2019}, and \cite{Chenetal2019}, but one should not expect to be so lucky in general.

Another approach uses \textit{synchronous couplings}: one diffusion process is initialized at $x$, and another process sharing the same Brownian motion is started at $x + \delta e^{(i)}$. The bound then depends on the coupling time of the two diffusions. This idea was exploited heavily in \cite{GorhMack2016} for instance to study derivative bounds for overdamped Langevin diffusions. There are other approaches aside from  synchronous couplings (we review them briefly in Section~\ref{sec:relatedwork}), each with its own merits and drawbacks. Ultimately, however, none are universally applicable to all problems, making derivative bounds a common bottleneck of the generator comparison approach. In this paper, we present a new way to bound  the left-hand side of \eqref{eq:main_difference}. Let us illustrate the main steps.

Fix a test function $h: \delta \Z^{d} \to \R$, defined only on the lattice $\delta \Z^{d}$ as opposed to $\R^{d}$ as before.    Now, instead of \eqref{eq:poisson_diffusion}, we consider the Poisson equation of the prelimit,
\begin{align}
G_X f_h(\delta k) = \E h(X) - h(\delta k), \quad k \in   \Z^{d}. \label{eq:poisson_CTMC}
\end{align}
Proposition 7.1 of \cite{Asmu2003} can be adapted to show that a solution  exists provided $\E \abs{h(X)} < \infty$. Furthermore, this solution   is unique up to a constant.  We are tempted to proceed analogously to \eqref{eq:main_difference} by taking expected values with respect to $Y$, but we cannot do so because $G_X f_h(\delta k)$ is not defined on $\R^{d} \setminus \delta \Z^{d}$. We get around this by interpolating the discrete Poisson equation. Namely, we introduce a spline $A$, which interpolates functions $f: \delta \Z^{d} \to \R$ and results in extended functions $A f: \R^{d} \to \R$.  By applying $A$ to both sides of \eqref{eq:poisson_CTMC}, we obtain the interpolated Poisson equation 
\begin{align*}
AG_X f_h(x)  =&\ \E h(X) - A h (x) , \quad x \in \R^{d}.
\end{align*}
Under some mild conditions on $A f_h(x)$, It\^{o}'s lemma implies $\E G_Y A f_h(Y) = 0$, and so we  take expected values with respect to $Y$ to arrive at
\begin{align}
\E h(X) - \E A h (Y) =&\ \E AG_X f_h(Y) \notag  \\
=&\   \E AG_X f_h(Y) - \E G_Y A f_h(Y).\label{eq:main_difference_discrete}
\end{align}
 To ensure that the convergence of \eqref{eq:main_difference_discrete} to zero implies the convergence of $X$ to $Y$, we again need to ensure that $h(\delta k)$ belongs to a rich-enough class of functions. We describe some convergence-determining classes of grid-restricted test functions in Section~\ref{sec:convdet}.    Lastly, to make the right-hand side of \eqref{eq:main_difference_discrete} comparable to \eqref{eq:main_difference}, we want to interchange $A$ and $G_X$.   This interchange is possible but results in some error; i.e.,\  $A G_X f_h(x) = G_X A f_h(x) + \text{ error}$.

 After this interchange, the right-hand side of \eqref{eq:main_difference_discrete} becomes analogous to \eqref{eq:main_difference} in the sense that the derivatives of $f_{h^{*}}(x)$ that appear in \eqref{eq:main_difference} are replaced by corresponding derivatives of $A f_h(x)$. Our choice of $A$ is such that the derivatives of $A f_h(x) $ correspond to finite differences of $f_h(\delta k)$, thus replacing the problem of establishing derivative bounds  by an analogous problem of bounding finite differences. The finite differences of $f_h(\delta k)$ are determined entirely by the Poisson equation, which itself is determined by the transition structure of the CTMC. As such, it is fitting to refer to the $k$th-order finite difference of $f_h(\delta k)$ as the \emph{$k$th-order Stein factor of the CTMC}.  We can bound these Stein factors by relying on the fact that
\begin{align}
f_h(\delta k) = \int_{0}^{\infty} \big( \E_{X(0) = \delta k} h(X(t)) - \E h(X) \big) dt, \quad k \in \Z^{d}, \label{eq:generalpoisson}
\end{align}
solves the Poisson equation  
and constructing synchronous couplings of the CTMC similar  to the diffusion synchronous couplings. Ways to verify that \eqref{eq:generalpoisson} is well defined are discussed in Section~\ref{sec:diffbounds}.

For ease of reference, we refer to our approach as the prelimit generator comparison approach, or simply \emph{prelimit approach}, and to the traditional approach based on \eqref{eq:poisson_diffusion} as the \emph{diffusion approach}. The prelimit and diffusion approaches are in some sense parallel approaches  with many conceptual similarities. If we choose $h^{*}(x)$ in \eqref{eq:main_difference} to equal $A h(x)$ from \eqref{eq:main_difference_discrete}, we see that the right-hand sides of \eqref{eq:main_difference} and \eqref{eq:main_difference_discrete} are equal. This means that any bound established via the diffusion approach should, in theory,   be attainable via the prelimit approach, and vice versa.

%For instance, a given CTMC can have multiple diffusion approximations. It was shown in \cite{BravDaiFang2020}  that using approximations with state-dependent diffusion coefficients, as opposed to constant ones, can reduce the approximation error by an order of magnitude. The diffusion approach requires establishing new derivative bounds for each different diffusion approximation, which is cumbersome. Furthermore, having a state-dependent diffusion coefficient complicates the diffusion synchronous coupling technique. Indeed, the results in \cite{Gorhetal2019} include diffusions with state-dependent diffusion coefficients, but do not include third derivative bounds, which are essential to the diffusion approach. In contrast,  the prelimit approach does not care about which diffusion approximation is used, and only requires us to establish difference bounds once. 
 In practice,   technical differences  can make the prelimit approach more attractive for some models.
 First, when working with models that have state-space collapse i.e.,\ when the dimension of the CTMC is higher than that of the diffusion the prelimit approach does not require one to bound the so-called $\E \abs{X_{\perp}}$, which is the distance between the stationary distribution of the CTMC and its projection onto the state-space collapse manifold. This is illustrated in more detail in \cite{Brav2022}, a companion paper in which the prelimit approach is applied to the join-the-shortest-queue model. Second, the diffusion approach can suffer from what we call ``misalignment of synchronous couplings,'' which can complicate the process of getting derivative bounds via diffusion synchronous couplings. We illustrate this issue  in Section~\ref{sec:compare_couplings} using a simple example.

%Looking back at \eqref{eq:first_deriv}, we see that using a diffusion synchronous coupling for first derivative bounds requires keeping track of two diffusion processes -- one initialized at $x$ and the other at $x + \varepsilon e^{(i)}$. Higher order derivatives require us to keep track of more processes, e.g.\ $\frac{d^2}{d x_i^2} f_{h^{*}}(x)$  is approximated by $\big( f_{h^{*}}(x + 2\varepsilon e^{(i)}) - 2 f_{h^{*}}(x + \varepsilon e^{(i)}) + f_{h^{*}}(x  ) \big) / \varepsilon^2 $ and  requires three processes.  First, while it is simple to construct a synchronous coupling when the diffusion process has a constant diffusion coefficient, it is not straightforward when the diffusion coefficient is state-dependent. Second, we give several examples of diffusion processes where diffusion synchronous couplings cannot be used to bound the third derivatives of $f_{h^{*}}(x)$, because the precise spacing of the coupled diffusions is violated over time.  We flesh out both of these points in more detail in Section~\ref{sec:drawbacks}.

Apart from showing how fast $X$ converges to $Y$, the prelimit Poisson equation  \eqref{eq:poisson_CTMC} can also be used to establish tightness of a family of steady-state distributions. Tightness has become an important property since the seminal work of \cite{GamaZeev2006}, which initiated a wave of research into justifying steady-state diffusion approximations of queueing systems; see, for instance, \cite{DaiDiekGao2014, BudhLee2009, ZhanZwar2008, Kats2010, YaoYe2012, Tezc2008, GamaStol2012}, and \cite{Gurv2014a}. Roughly speaking, process-level convergence of the CTMC to a diffusion  combined with tightness of the CTMC stationary distributions enables one to perform a limit-interchange argument to conclude convergence of steady-state distributions. The bottleneck is usually proving tightness, which  has become synonymous with steady-state convergence.

We can use \eqref{eq:poisson_CTMC} to prove tightness as follows. Let  $x_{\infty}$ be the fluid equilibrium of the CTMC, and assume for simplicity that $x_{\infty} \in \delta \Z^{d}$ (if not, consider the nearest point in $\delta \Z^{d}$). Pick $h(\delta k) = \abs{\delta k - x_{\infty}}$ and evaluate the Poisson equation at the point $\delta k = x_{\infty}$ to get 
\begin{align*}
\E \abs{X - x_{\infty} } = G_{X} f_h(x_{\infty}).
\end{align*}
The right-hand side typically contains CTMC Stein factors up to the second order.  We give an example in Section~\ref{sec:tight}. Proving tightness is therefore equivalent to bounding these factors at the \emph{single point} $x_{\infty}$. In contrast, bounding the approximation error of $Y$ requires third-order Stein factor bounds on the \emph{entire support} of $Y$. This  highlights the extra work needed for convergence rates as opposed to convergence alone.

%To showcase the discrete Poisson approach, we use the join-the-shortest queue (JSQ) system as a working example in this paper. Apart from the renewed attention it received in the recent years  \cite{GamaEsch2018, Brav2019,BaneMukh2019, BaneMukh2019b, GuptWalt2019, LiuYing2019b}, the JSQ model is significant from a technical standpoint because it is a non-trivial multi-dimensional model with a state-space collapse component. This makes it a good example to show off the power of the discrete approach.

The idea of interpolating the discrete Poisson equation can be applied more broadly to the problem of comparing discrete and continuous distributions using Stein's method. To the author's knowledge, anytime Stein's method has been invoked for a discrete-versus-continuous random variable comparison, the starting point has always been the differential equation for the continuous random variable. Furthermore, in most applications of the method,  the starting point has been the Stein/Poisson equation for the limiting distribution, whereas we start with the prelimit. 

To summarize, our main contribution is the prelimit approach, which depends on two technical components. First, we establish the existence of an interpolator $A$ that satisfies certain convenient properties. Theorem~\ref{thm:1dinterpolant_def} contains the one-dimensional result, which is generalized to multiple dimensions in  Theorem~\ref{thm:interpolant_def}. Second, we  describe  the error of interchanging $A$ with $G_{X}$.   Proposition~\ref{lem:interchanged1}  contains  the one-dimensional result, while Proposition~\ref{lem:mdiminterchange} is the  multidimensional generalization.

After illustrating the general framework, we apply it to the $M/M/1$ queueing system to showcase the prelimit approach. The steady-state customer count in the $M/M/1$ model is geometrically distributed, and is approximated by the exponential distribution; convergence rates  are presented in  Theorem~\ref{thm:mm1}.  The  $M/M/1$ system is chosen purely for illustrative purposes because of its simplicity, and our convergence rates actually have   alternative derivations. For example, one can use  existing results on Stein's method for the exponential distribution in  Theorem 3.1 of \cite{PekoRoll2011} or Theorem 5.11 of \cite{Ross2011}. Furthermore, the Poisson equation for the $M/M/1$ system is the same as the Stein equation for the geometric distribution, which was first obtained in \cite{peko1996} and is also a special case of the Pascal Stein equation considered in \cite{Scho2001}. Bounds on  geometric Stein factors have also been obtained in \cite{Daly2008}. We compare existing bounds with our own in Section~\ref{sec:diffbounds}.

It is important to add that using CTMC synchronous couplings dates back to \cite{Barb1988}, which was  the first paper to connect Stein's method to Markov chains (the author of that paper did not use the language ``synchronous coupling''). In that work, the author viewed the Poisson distribution as the steady-state distribution of the infinite server queue. Later, the application of Stein's method to birth-death processes received a thorough treatment in \cite{BrowXia2001}.  A more recent example of using CTMC synchronous couplings can be found in \cite{BarbLuczXia2018a,BarbLuczXia2018b}.

The remainder of the paper is structured as follows. In Section~\ref{sec:prelimapproach} we introduce the technical components of the prelimit approach. We then apply the prelimit approach to the $M/M/1$ model and illustrate the synchronous coupling idea in Section~\ref{sec:diffbounds}. We discuss the issue of misalignment of synchronous couplings in Section~\ref{sec:compare_couplings} and conclude in Section~\ref{sec:conclusion}.

\subsection{Related Work on Derivative Bounds}
\label{sec:relatedwork}
Let us briefly discuss several recent works on ways to obtain derivative bounds.   In \cite{GorhMack2016} the authors used synchronous couplings to study derivative bounds for overdamped Langevin diffusions with strongly concave drifts. Later in  \cite{Gorhetal2019}, the authors relaxed the strongly concave drift assumption to a dissipativity condition and used a combination of synchronous couplings and reflection couplings studied in \cite{Eber2016} and \cite{Wang2016} to establish derivative bounds for a class of fast-coupling diffusions. In \cite{ErdoMackSham2019} the authors establish derivative bounds for an even larger class of diffusions, but still require a dissipativity condition. The strong concavity and dissipativity conditions both imply that the diffusion generator $G_Y$ satisfies 
\begin{align}
G_Y V(x) \leq - \alpha \abs{x}^2 + \beta, \quad x \in \R^{d}, \label{eq:dissip}
\end{align} 
where  $\abs{x}$ denotes the Euclidean norm, $V(x) = \abs{x}^2$, and $\alpha, \beta$ are some positive constants. Condition \eqref{eq:dissip} is also known as $V$-exponential ergodicity (with $V(x) = \abs{x}^2$), see \cite{MeynTwee1993b}.

 While the aforementioned papers contain a large list of applications, their results are not directly applicable to many queueing settings because \eqref{eq:dissip} does not hold there. Even one of the most basic diffusion processes in queueing, the piecewise  Ohrnstein-Uhlenbeck process used for approximating the many-server queue  in \cite{BravDaiFeng2016}, does not satisfy \eqref{eq:dissip}.
% The derivative bounds of \cite{GorhMack2016,Gorhetal2019} depend on the magnitude of the derivatives of the drift. This means that even if the drift in \cite{BravDaiFeng2016} is smoothed  and made differentiable, the bounds in \cite{GorhMack2016,Gorhetal2019} would not be tight enough and could not be used to recover the convergence rate result of \cite{BravDaiFeng2016}. 
 Furthermore, the results in these papers hold only  for diffusions on the entire space $\R^{d}$. This excludes diffusions with reflecting boundary conditions, such as reflecting Brownian motions that appear as heavy-traffic limits for networks of single-server queueing systems.

Another approach to getting derivative bounds was proposed in \cite{Gurv2014}, where the author used a priori Schauder estimates from PDE theory to bound the derivatives of $f_{h^{*}}(x)$ in terms of $f_{h^{*}}(x)$ and $h(x)$. He then bounded $f_{h^{*}}(x)$ by a Lyapunov function satisfying an exponential ergodicity condition for the diffusion. This approach requires finding a Lyapunov function  satisfying an exponential ergodicity condition, which   typically requires significant effort, e.g.\ \cite{DiekGao2013, Gurv2014}. Furthermore, in the case of a diffusion with a reflecting boundary, the complexity of the PDE machinery used makes it nontrivial to  trace how the a priori Schauder estimates depend on the primitives of the diffusion process.  

Most recently, another approach to getting derivative bounds based on Bismut's formula from Malliavin calculus was  proposed in \cite{FangShaoXu2018}.  The authors required the diffusion coefficient to be constant, and the assumptions imposed on the drift were similar to those  in \cite{GorhMack2016}.

\subsection{Notation}
For any $B \subset \R^{d}$,  let $\text{Conv}(B)$ denote its convex hull.  We use $\Z$ to denote the set of integers  and let $\N = \{0,1,2,\ldots\}$. For any   $k \in \N $ and  $B \subset \R^{d}$, we let $C^{k}(B)$ be the set of all $k$-times continuously differentiable functions $f: B \to \R$.   Given a stochastic process $\{Z(t)\} \in D$ and a functional $f: D \to \R$, we write $\E_{x}(f(Z))$ to denote $\E(f(Z)\ |\ Z(0) = x)$. We let $e \in \R^{d}$ be the vector whose elements all equal  $1$  and let $e^{(i)}$ be the element with $1$ in the $i$th entry  and zeros otherwise. For any $\delta > 0$ and integer $d > 0$, we let $\delta \Z^{d} = \{\delta k :\ k \in \Z^{d}\}$ and define $\delta \N^{d}$ similarly. For any function $f: \delta \Z^{d} \to \R$, we define the forward difference operator in the $i$th direction as
\begin{align*}
\Delta_{i} f(\delta k) = f \big( \delta (k+e^{(i)}) \big) - f(\delta k), \quad k \in \Z^{d}, \  1 \leq i \leq d,
\end{align*}
and for $j \geq 0$, we define 
\begin{align}
\Delta_i^{j+1} f(\delta k) = \Delta_{i}^{j} f(\delta(k+e^{(i)})) - \Delta_{i}^{j} f(\delta k), \label{eq:diffdef}
\end{align}
with the convention that $\Delta_i^{0} f(\delta k) = f(\delta k)$. For a vector $a \in \N^{d}$, we also let 
\begin{align*}
\Delta^{a}  f(\delta k) =&\ \Delta_{1}^{a_1} \ldots \Delta_{d}^{a_d} f(\delta k),
\end{align*} 
and if  $f: \R^{d} \to \R$, then
\begin{align*}
\frac{\partial^{a}}{\partial x^{a}} f(x) =&\ \frac{\partial^{a_1}}{\partial x_1^{a_1}} \ldots \frac{\partial^{a_d}}{\partial x_d^{a_d}} f(x),
\end{align*}
and we adopt the convention that $\frac{\partial^{0}}{\partial x^{0}} f(x) = f(x)$. For any $x \in \R^{d}$, we define $\norm{x}_{1} = \sum_{i=1}^{d} \abs{x_i}$  and write $\abs{x}$ to denote the Euclidean norm. Throughout the paper we will often use $C$ to denote a generic positive constant that may change from line to line  and that will    be independent of any parameters not explicitly specified.  For a random variable $X$, we write $\supp(X)$ to denote the support of $X$.

\section{The Prelimit Generator Comparison Approach}
\label{sec:prelimapproach}
In this section, we work out the technical details of the prelimit approach. We begin by   introducing the interpolation operator $A$ in Section~\ref{sec:interpol}. We follow this with a discussion of convergence-determining classes in Section~\ref{sec:convdet}.  Then, we write the form of $A G_{X} f(x)$ in a manner that easily lends itself to analysis. Informally, we refer to this  as interchanging $A$ with $G_{X}$. Bounded and unbounded domains require separate consideration. We treat unbounded domains in Section~\ref{sec:noref} and treat one example of a bounded domain in Section~\ref{sec:bounded}.   To minimize notational burden, we restrict our discussion to one-dimensional CTMCs. In multiple dimensions, the results are analogous from a technical perspective, but may be harder to parse at first read. We therefore postpone the multidimensional discussion to the appendix, in which multidimensional interpolation is discussed in Appendix~\ref{sec:interpolation}, and multidimensional interchange is left to Appendix~\ref{sec:mdimunbound}.

\subsection{The Interpolator}
\label{sec:interpol}
The objective of this section is to state Theorem~\ref{thm:1dinterpolant_def}.  Fix $\delta > 0$, and for $x \in \R$ define $k(x) = \lfloor x/\delta\rfloor$. Let $K \subset \R$ be a possibly unbounded interval and define 
\begin{align*}
K_{4} = \{x \in K \cap \delta \Z : (x + 4\delta) \in K \cap \delta \Z\}.
\end{align*}
For example, if $K = (-\infty, \infty)$, then $K_4 = \delta \Z$. Let $f: K \cap \delta \Z \to \R $ be the function we want to extend to the continuum. 
%One may be familiar with the cubic Hermite spline,  which can certainly extend the function. However, we need a higher order spline  because 
We interpolate the function using splines, which  are standard tools in numerical analysis; see for instance Section 8 in  \cite{Kres1998}. Instead of the popular cubic spline, which only results in a $C^1$ interpolant, we craft a  degree-7 spline so that our extension is thrice continuously differentiable.   Define 
\begin{align*}
A f(x) =&\ P_{k(x)}(x), \quad \text{ where }  \quad P_{k}(x)   =   \sum_{i=0}^{4} \alpha^{k}_{k+i}(x) f(\delta (k+i)), \quad x \in \R.
\end{align*}
Each $P_{k}(x)$ is a degree-7 polynomial and is best understood as a weighted sum of $f(\delta k), \ldots, f(\delta (k+4))$  with weights $\alpha^{k}_{k}(x), \ldots, \alpha^{k}_{k+4}(x)$.  The precise form of $P_{k}(x)$ is distracting,  so we state it in Appendix~\ref{sec:interpolation}. The following result summarizes the key properties  we require of $A f(x)$ and the weights $\alpha^{k}_{k+i}(x)$. 
\begin{theorem} \label{thm:1dinterpolant_def}
Given $f: K \cap \delta \Z \to \R$, the function   
\begin{align}
Af(x) = \sum_{i=0}^{4} \alpha^{k(x)}_{k(x)+i}(x)f(\delta(k(x)+i)), \quad x \in \text{Conv}(K_4), \label{eq:af}
\end{align}
belongs to  $ C^{3}(\text{Conv}(K_4))$ and is infinitely differentiable on $\text{Conv}(K_4) \setminus K_4$. Furthermore, 
\begin{align}
A f(\delta k) = f(\delta k), \quad \delta k \in K_4, \label{eq:interpolates}
\end{align}
and the derivatives of $A f(x)$ are bounded by the corresponding finite differences of $f(\delta k)$. Namely, there exists $C > 0$ independent of $x$, $f(\cdot)$ and $\delta$  such that 
\begin{align}
\Big|  \frac{\partial^{a}}{\partial x^{a}} Af(x) \Big|   \leq&\  C \delta^{-a}  \max_{\substack{ 0 \leq i  \leq 4-a  }} \abs{\Delta^{a} f(\delta (k(x)+i))}, \quad x \in \text{Conv}(K_4),\ 0 \leq a \leq 3, \label{eq:multibound}
\end{align} 
and \eqref{eq:multibound} also holds  for $x \in  \text{Conv}(K_4) \setminus K_4$ when $a = 4$. Additionally, the weights  $\big\{\alpha_{k+i}^{k}: \R \to \R \ | \ k \in \Z,\ i = 0,1,2,3,4\big\}$ 
are degree-$7$ polynomials in $(x-\delta k)/ \delta$ whose coefficients do not depend on $k$ or $\delta$. They satisfy  
\begin{align}
&\alpha_{k}^{k}(\delta k) = 1, \quad \text{ and } \quad \alpha_{k+i}^{k} (\delta k) = 0, \quad &k \in \Z,\ i = 1,2,3,4, \label{eq:alphas_interpolate} \\
&\sum_{i=0}^{4} \alpha^{k}_{k+i}(x) = 1, \quad &k \in \Z,\ x \in \R, \label{eq:weights_sum_one}
\end{align}
and also the following translational invariance property:
\begin{align}
\alpha^{k+j}_{k+j+i}(x+ \delta j)  = \alpha^{k}_{k+i}(x),\quad i,j,k \in \Z, \ x \in \R. \label{eq:weights}
\end{align}
\end{theorem}
 Theorem~\ref{thm:1dinterpolant_def} is proved in Appendix~\ref{sec:interpolation} and follows directly from the form of $P_{k}(x)$ stated there. From \eqref{eq:multibound} we see that the reason $P_{k}(x)$ depends on  $f(\delta k), \ldots, f(\delta (k+4))$, as opposed to  also depending on $f(\delta (k+5))$, is that we want $\frac{\partial^{a}}{\partial x^{a}} Af(x)$ to be related to $\Delta^{a} f(\delta k(x))$ for $0 \leq a \leq 4$, and we do not care what happens beyond the fourth derivative. In theory, one can  make $A f(x)$ as differentiable as is needed by using a higher degree polynomial $P_{k}(x)$.

\subsection{Convergence-Determining Classes}
\label{sec:convdet}
We mentioned in the introduction that when one uses the diffusion approach,  $\lipone$ is a commonly used convergence-determining class. In this section we discuss two convergence-determining classes of grid-valued functions that can be used with the prelimit approach. Lemma~\ref{lem:convdet} below presents the main result of this section.
 
Recall our convention of using a star superscript to emphasize that a function is defined on the  continuum.  
Given two random variables $U,V \in \R^{d}$ and a class of functions $\mathcal{H} = \{h^{*}: \R^{d} \to \R\}$, we define 
\begin{align*}
d_{\mathcal{H}}(U,V) = \sup_{h^{*} \in \mathcal{H}} \Big| \E h^*(U) - \E h^*(V) \Big|.
\end{align*}
We already said that $\lipone$ is a convergence-determining class because $d_{\lipone}(U,V) \to 0$ implies $U$ converges to $V$ in distribution. There are, of course, other convergence-determining classes. For instance, it was shown in Lemma 2.2 of \cite{GorhMack2016} that if
\begin{align*}
\mathcal{H} = \mathcal{M} = \Big\{h^*: \R^{d}  \to \R : \Big| \frac{\partial^{a}}{\partial x^a} h^{*}(x)  \Big| \leq 1, \quad 1 \leq \norm{a}_{1} \leq 3  \Big\},
\end{align*}  
then $d_{\mathcal{M}}(U,V) \to 0$ also implies convergence in distribution.  

Both $\lipone$ and $\mathcal{M}$ are classes of functions defined on $\R^{d}$, but  the prelimit approach works with functions defined only on $\delta \Z^{d}$. To mimic the two classes, we define 
\begin{align*}
\text{dLip}(1) =&\ \{h: \delta \Z^{d} \to \R : \abs{\Delta_{j} h(\delta k)} \leq \delta,\ 1 \leq j \leq d,\  k \in \delta \Z^{d} \},\\
\mathcal{M}_{disc}(C)   =&\ \{h: \delta \Z^{d} \to \R : \abs{\Delta^{a} h(\delta k)} \leq C \delta^{\norm{a}_{1}} ,\ 1 \leq \norm{a}_{1} \leq 3,\  k \in \delta \Z^{d} \}.
\end{align*}
The following lemma relates $d_{\lipone}(U,V)$ and $d_{\mathcal{M}}(U,V)$ to their grid-restricted counterparts. The lemma involves the multidimensional interpolator, which we have not yet formally introduced. However, that  does not preclude an understanding of the lemma, which is proved   in Section~\ref{sec:proofconvdet}.
\begin{lemma}
\label{lem:convdet}
Let $U \in \delta \Z^{d}$ and $V \in \R^{d}$ be two random vectors. For any $h^{*}: \R^{d} \to \R$, let   $h: \delta \Z^{d} \to \R$ be the restriction of $h^{*}(x)$ to $\delta \Z^{d}$. Let $A h: \R^{d} \to \R$ be defined by \eqref{eq:af} when $d =1$, and by \eqref{eq:af2} when $d > 1$. Then there exists a constant $C > 0$ such that 
\begin{align*}
\abs{\E h^{*}(U) - \E h^{*}(V)} \leq \abs{\E h(U) - \E A h(V)} + C \delta \max_{\substack{  1 \leq j \leq d   }} \sup_{\substack{  x \in \R^{d} }} \bigg| \frac{\partial}{\partial x_{j}} h^{*}(x)  \bigg|.
\end{align*}
As a consequence, there exists a constant $C'>0$ such that 
\begin{align*}
d_{ \lipone }(U,V) \leq \sup_{h \in \dlipone} \abs{\E h(U) - \E A h(V)} + C \delta,\\ 
d_{\mathcal{M}}(U,V) \leq \sup_{h \in \mathcal{M}_{disc}(C' )} \abs{\E h(U) - \E A h(V)} + C \delta.
\end{align*} 
\end{lemma}
 
\subsection{Interchange for Unbounded Domains}
\label{sec:noref}
 Assume $G_{X} f(\delta k)$ is defined for all $k \in \Z$; i.e.,\ the CTMC lives on $\delta \Z$. The interchange result for $A$ and $G_X$ is given in Proposition~\ref{lem:interchanged1} below. We then apply this result to characterize the approximation error between $X$ and its diffusion approximation $Y$ in \eqref{eq:errordiff}.

  Define $\beta_{\ell}(\delta k) = q_{\delta k,\delta (k+\ell)}$ for $k,\ell \in \Z$, where $q_{\delta k, \delta k'}$ are the CTMC transition rates.
%\begin{align*}
%\beta_{\ell}(\delta k) = q_{\delta k,\delta (k+\ell)}, \quad k, \ell \in \Z.
%\end{align*}
 Then
\begin{align*}
G_X f(\delta k) =  \sum_{k' \in \Z } q_{\delta k,\delta k'} (f(\delta k')-f(\delta k)) =&\ \sum_{\ell \in \Z} \beta_{\ell}(\delta k) (f(\delta (k+\ell))-f(\delta k)), \quad k \in \Z.
\end{align*}
Fix $h : \delta \Z \to \R$ with $\E \abs{h(X)} < \infty$. Since $A$ is a linear operator, we apply to both sides of the CTMC Poisson equation to get  
\begin{align*}
 A G_{X} f_h(x) = A(\E h(X) - h) (x) = \E A h(X) - A h(x). 
\end{align*}
The following result says $A G_{X} f(x) = G_{X} A f(x) + \text{error}(x)$ and characterizes the error term. We prove it in Section~\ref{sec:mdimunbound}  by proving the multidimensional version, Proposition~\ref{lem:mdiminterchange}, there.
\begin{proposition} \label{lem:interchanged1}
Fix $f: \delta \Z \to \R$ and assume that  $G_{X} f(\delta k)$ is defined on all of $\delta \Z$.  Assume also that 
\begin{align}
 \sum_{\ell \in \Z} \abs{\beta_{\ell}(\delta k) (f(\delta (k+\ell))-f(\delta k))} < \infty, \quad k \in \Z, \label{eq:intergrab}
\end{align}
which is trivially satisfied  when the number of transitions from each state is finite.  Then 
\begin{align}
 A G_{X} f(x) =&\  \sum_{\ell \in \Z}  A \beta_{\ell}(x)  \big( A f(x+\delta \ell) - A f(x)\big) + \varepsilon(x), \quad x \in \R.    \label{eq:intererror}
\end{align}
 The error $\varepsilon(x)$  satisfies 
\begin{align}
\varepsilon(x) =&\ \sum_{\ell \in \Z}  \sum_{i=0}^{4} \alpha^{k(x)}_{k(x)+i}(x)\Big(\beta_{\ell}\big(\delta (k(x)+i)\big) - A \beta_{\ell}(x)\Big) \notag \\
 &   \times  \Big(1(\ell > 0) \sum_{j=0}^{i-1} \sum_{m=0}^{\ell-1} \Delta^{2} f\big( \delta(k(x)+ m  + i )\big) - 1(\ell < 0) \sum_{j=0}^{i-1} \sum_{m=\ell}^{-1} \Delta^{2} f\big( \delta(k(x)+ m  + i )\big) \Big). \label{eq:vareps}
\end{align}
%\begin{align*}
%\varepsilon(x) =&\ \sum_{\ell \in \Z}  \sum_{i=0}^{4} \alpha^{k(x)}_{k(x)+i}(x)\Big(\beta_{\ell}\big(\delta (k(x)+i)\big) - A \beta_{\ell}(x)\Big) \notag  \\
%& \hspace{1.5cm} \times \Big(f\big(\delta (k(x)+\ell+i)\big)-f\big(\delta (k(x)+i)\big) - \big(f\big(\delta (k(x)+\ell )\big)-f\big(\delta k(x) \big)\big)\Big).
%\end{align*}
\end{proposition}
Let us now fix a CTMC and derive a diffusion approximation $Y$ for it. We also characterize the approximation error. Fix $h(\delta k)$ and consider   \eqref{eq:intererror} with $f(\delta k) = f_h(\delta k)$ (we assume \eqref{eq:intergrab} holds).  First, we  apply Taylor expansion to  $\big( A f_{h}(x+\delta \ell) - A f_{h}(x)\big)$ to get
\begin{align*}
A G_{X} f_h(x) =&\ (A f_h)'(x) \delta \sum_{\ell \in \Z} \ell A \beta_{\ell}(x)  + \frac{1}{2} (A f_h)''(x) \delta^{2} \sum_{\ell \in \Z} \ell^{2} A \beta_{\ell}(x) \\
&+ \frac{1}{6} \delta^{3} \sum_{\ell \in \Z} \ell^{3} A \beta_{\ell}(x) (A f_h)'''(\xi_{\ell}(x)) + \varepsilon(x), 
\end{align*} 
where $\xi_{\ell}(x)$ is  between $x$ and $x + \delta \ell$. To approximate $X$, we set
\begin{align*}
b(x) = \delta \sum_{\ell \in \Z} \ell A \beta_{\ell}(x), \quad \text{ and } \quad a(x) = \delta^2 \sum_{\ell \in \Z} \ell^2 A \beta_{\ell}(x), \quad x \in \R,
\end{align*}
and consider the diffusion process
\begin{align}
Y(t) = Y(0) + \int_{0}^{t} b(Y(s)) ds + \int_{0}^{t} \sqrt{a(Y(s))} d W(s), \label{eq:sdediffusion}
\end{align}
where $\{W(t)\}$ is standard Brownian motion. The generator of this diffusion is  $G_{Y} f(x) = b(x) f'(x) + \frac{1}{2} a(x) f''(x)$ and its stationary distribution has density  $\frac{\kappa}{a(x)} \exp\big(\int_{0}^{x} \frac{2 b(y)}{a(y)} dy\big)$, where $\kappa$  is a normalizing constant that we assume to be finite. Let $Y$ be the random variable having this density.  It\^{o}'s lemma  tells us that for any $f \in C^{2}(\R)$, 
\begin{align*}
 \E_{x} f(Y(t)) - \E_{x} f(Y(0)) =&\ \E_{x} \int_{0}^{t} G_Y f(Y(s)) ds, \quad t > 0.
\end{align*}
Provided $\E |f(Y)|<\infty$,  we can initialize $Y(0) \stackrel{d}{=} Y $ to get
\begin{align*}
\E \Big[  \int_{0}^{t} G_Y f(Y(s)) ds  \Big| Y(0) \stackrel{d}{=} Y \Big] = 0.
\end{align*}
 If we further assume that $\E \abs{G_Y f(Y)} < \infty$, then we can apply the Fubini-Tonelli theorem to interchange the integral and expectation above and conclude that 
\begin{align}
\E G_{Y} f(Y) = 0. \label{eq:ito}
\end{align}
%random variable $Y$ on $\R$ with density  $\frac{\kappa}{a(x)} \exp\big(\int_{0}^{x} \frac{2 b(y)}{a(y)} dy\big)$, where $\kappa$  is a normalizing constant that we assume to be finite. Integration by parts yields 
%\begin{align}
%\E b(Y) f(Y) + \frac{1}{2} \E a(Y) f'(Y) = 0 \label{eq:gy}
%\end{align}
%for any $f(x)$ for which the expectations above exist and for which $\lim_{x \to \pm \infty} f(x) \exp\big(\int_{0}^{x} \frac{2 b(y)}{a(y)} dy\big) = 0$.  Another way to view $Y$ is as the stationary distribution of the  diffusion process  
 Now, provided that \eqref{eq:ito} holds  with $A f_h(x)$ in place of $f(x)$ there, we get
\begin{align}
 \E A h(X) - \E A h(Y) =&\ \E A G_{X} f_h(Y) - \E G_{Y} A f_{h}(Y) \notag \\
 =&\   \frac{1}{6} \delta^{3} \E \sum_{\ell \in \Z} \ell^{3} A \beta_{\ell}(Y) (A f_h)'''(\xi_{\ell}(Y)) + \E \varepsilon(Y). \label{eq:errordiff}
\end{align}
 The bounds on  $(A f_{h})'''(x)$ from Theorem~\ref{thm:1dinterpolant_def} imply 
\begin{align*}
\frac{1}{6} \delta^{3} \Big| \E \sum_{\ell \in \Z} \ell^{3} A \beta_{\ell}(Y) (A f_h)'''(\xi_{\ell}(Y)) \Big| \leq&\ C   \Big| \E \sum_{\ell \in \Z} \ell^{3} A \beta_{\ell}(Y)  \max_{\substack{ 0 \leq i  \leq 1  }} \abs{\Delta^{3} f_h(\delta (k(\xi_{\ell}(Y))+i))} \Big|.
\end{align*}
In other words, the term above depends on $\ell^{3} A \beta_{\ell}(Y)$ and third-order Stein factors. The second  term in \eqref{eq:errordiff} is $\E \varepsilon(Y)$. We recall $\varepsilon(x)$ below for convenience:
\begin{align*}
\varepsilon(x) =&\ \sum_{\ell \in \Z}  \sum_{i=0}^{4} \alpha^{k(x)}_{k(x)+i}(x)\Big(\beta_{\ell}\big(\delta (k(x)+i)\big) - A \beta_{\ell}(x)\Big) \notag \\
 &   \times  \Big(1(\ell > 0) \sum_{j=0}^{i-1} \sum_{m=0}^{\ell-1} \Delta^{2} f\big( \delta(k(x)+ m  + i )\big) - 1(\ell < 0) \sum_{j=0}^{i-1} \sum_{m=\ell}^{-1} \Delta^{2} f\big( \delta(k(x)+ m  + i )\big) \Big).
\end{align*}
First,  the fact that $\alpha^{k }_{k +i}(x)$ is a polynomial in $(x-\delta k)/\delta$ implies  $\sup_{x \in \R} \big|\alpha^{k(x)}_{k(x)+i}(x)\big|$ is bounded by a constant independent of $\delta, x$,  or any other parameters.   Second,  the fact that $0 \leq i \leq 4$ and the mean value theorem imply that 
\begin{align*}
\abs{\beta_{\ell}\big( \delta(k(x)  + i)\big) - A \beta_{\ell}(x)} = \abs{A \beta_{\ell}\big( \delta(k(x)  + i)\big) - A \beta_{\ell}(x)} \leq  4\delta (A \beta_{\ell})'(\xi'_{\ell}(x)).
\end{align*}
 Therefore, provided that the transition rates $\beta_{\ell}(\cdot)$ of the CTMC do not vary too much, e.g.,\ they are Lipschitz, the term above can be controlled, so bounding \eqref{eq:errordiff} comes down to bounding $\Delta^{2} f_h(\delta k)$ and $\Delta^{3} f_h(\delta k)$. 
% Practice has shown that higher-order differences of $f_h(\delta k)$ are generally smaller. For example, in  Proposition~\ref{lem:mm1dif} of Section~\ref{sec:diffbounds} we show that for the single-server queueing system,  the ratio $\E \Delta^{j+1} f_h(X)/\E \Delta^{j} f_h(X)$ is on the order of $(1-\rho)$, where $\rho$ is the system utilization.

\subsection{Interchange for a Bounded Domain}
\label{sec:bounded}
When the domain of the CTMC is bounded, Proposition~\ref{lem:interchanged1} must be modified slightly to account for the boundary of the domain. In this section we illustrate this using the example of the birth-death process defined by the generator
\begin{align}
G_{X} f(\delta k) = \lambda \Delta f(\delta k) - \mu 1(k > 0)  \Delta f(\delta (k-1)), \quad k \in \N. \label{eq:mm1gen}
\end{align}
 This generator corresponds to the customer count, scaled by $\delta$, in a single-server queue where customers arrive according to a  Poisson process with rate $\lambda$ and service  times are exponentially distributed with rate $\mu$. Such a system is also known as the $M/M/1$ queueing system. The quantity $\rho  = \lambda/\mu$ is the system utilization. In steady state, the customer count is geometrically distributed provided that $\rho < 1$. It is also well known that as $\rho \to 1$, the customer count can be approximated by an exponential random variable.  A recent application of Stein's method in \cite{GaunWalt2020} establishes convergence rates of the waiting time distribution in the $M/G/1$ system (which is more general than the $M/M/1$ system) to the exponential distribution. Another example of a CTMC with a bounded domain can be found in \cite{Brav2022}. 

The $M/M/1$ system is restricted to the non-negative integers.   Let us see how this boundary affects the interchange of $A$ and $G_X$.  Fix $f: \delta \N \to \R$ and consider $A G_X f(x)$, which is defined for  $x \in [0,\infty)$. For $x  \geq \delta$, the proof of Proposition~\ref{lem:interchanged1} can be repeated to see that 
 \begin{align*} 
 A G_{X} f(x) =   \lambda  \big(A f\big(x + \delta\big)-A f(x)\big) + \mu \big(A f\big(x - \delta\big)-A f(x)\big) + \varepsilon(x).  
 \end{align*}
 In fact, $\varepsilon(x) = 0$ for $x \geq \delta$ because the birth and death rates $\lambda$ and $\mu$ are constant. However, the equality above does not hold when $x \in [0,\delta)$ because $A f\big(x - \delta\big)$ is not defined there. To see this, we recall that $A f(x) = \sum_{i=0}^{4} \alpha^{k(x)}_{k(x)+i}(x)f(\delta(k(x)+i)) $, and that $f(\delta k)$ is only defined for $k \geq 0$. Our restriction of $f(\delta k)$ to  $k \geq 0$ is not artificial  because instead of $f(\delta k)$ we intend to use  the $M/M/1$  Poisson equation solution, which is defined only on $\delta \N$. 
 
To resolve this, we extend the definition $f(\delta k)$ to $k = -1$. To motivate the extension, fix $x \in [0,\delta)$ and consider 
\begin{align}
 A G_{X} f(x) =  \sum_{i=0}^{4} \alpha^{0}_{i}(x)\lambda \Big( f\big(\delta(i+1)\big) - f(i)\Big) + \sum_{i=1}^{4} \alpha^{0}_{i}(x)\mu \Big( f\big(\delta(i-1)\big) - f(i)\Big). \label{eq:mm1bdry}
\end{align}
Using the translational invariance property presented in \eqref{eq:weights} of  Theorem~\ref{thm:1dinterpolant_def}, it follows that 
\begin{align*}
\sum_{i=0}^{4} \alpha^{0}_{i}(x)\lambda \Big( f\big(\delta(i+1)\big) - f(i)\Big) =&\ \lambda \Big( \sum_{i=0}^{4} \alpha^{1}_{1+i}(x+\delta )f\big(\delta(i+1)\big) - \sum_{i=0}^{4} \alpha^{0}_{i}(x)f(i)\Big) \\
=&\ \lambda\big( A f(x+\delta) - A f(x)\big).
\end{align*}
We wish to do the same thing to the second term in \eqref{eq:mm1bdry}, but we cannot because the summation there starts from $i = 1$, not $i = 0$. Note that 
\begin{align*}
\sum_{i=1}^{4} \alpha^{0}_{i}(x)\mu \Big( f\big(\delta(i-1)\big) - f(i)\Big) =&\ \alpha^{0}_{0}(x) \mu \big(f(0) - f(0)\big) + \sum_{i=1}^{4} \alpha^{0}_{i}(x)\mu \Big( f\big(\delta(i-1)\big) - f(i)\Big)\\
=&\ \mu  \Big( \alpha^{0}_{0}(x)f(0) +    \sum_{i=1}^{4} \alpha^{0}_{i}(x) f\big(\delta(i-1)\big) \Big)  - \mu A f(x),
\end{align*}
which motivates us to define
\begin{align*}
\widehat f(\delta k) = 
\begin{cases}
f(\delta k), \quad k \geq 0, \\
f(0), \quad k = -1.
\end{cases}
\end{align*} 
It follows that $\mu  \Big( \alpha^{0}_{0}(x)f(0) +    \sum_{i=1}^{4} \alpha^{0}_{i}(x) f\big(\delta(i-1)\big) \Big) $ equals
\begin{align*}
   \mu \sum_{i=0}^{4} \alpha^{0}_{i}(x)  \widehat f\big(\delta(i-1)\big)   
=  \mu \sum_{i=0}^{4} \alpha^{-1}_{-1+ i}(x - \delta)  \widehat f\big(\delta(i-1)\big)  =  \mu  A \widehat f(x - \delta)   .
\end{align*}
Since $A f(x) = A \widehat f(x)$ for $x \geq 0$, we conclude that 
\begin{align}
 A G_{X} f(x) = \lambda  \big( A \widehat f(x + \delta) - A \widehat f(x)\big) + \mu \big( A \widehat f(x - \delta) - A \widehat f(x)\big), \quad x \geq 0. \label{eq:interchangemm1}
\end{align}
This result resembles  Proposition~\ref{lem:interchanged1}, but uses the extension $\widehat f(\delta k)$ instead of $f(\delta k)$. We comment more on our choice of extension at the end of this section.  To derive the diffusion approximation, we perform Taylor expansion. For  $x \geq 0$, 
\begin{align}
A G_{X} f_h(x) =&\ \lambda  ( A  \widehat f_h(x+\delta ) - A   \widehat f_h(x))  + \mu  ( A  \widehat f_h(x-\delta ) - A   \widehat f_h(x))   \notag \\
=&\ \delta (\lambda - \mu) (A  f_h)'(x) + \frac{1}{2}\delta^{2} (\lambda + \mu) (A  f_h)''(x) \notag \\
&+ \frac{1}{6} \delta^{3} \big( \lambda (A  f_h)'''(\xi_{1}(x)) - \mu (A  \widehat f_h)'''(\xi_{-1}(x))\big).   \label{eq:rbmgx}
\end{align}
In the second equality, we used the fact that $A  \widehat f_h(x) = A  f_h(x)$ for $x \geq 0$.  The diffusion approximation is driven by the first- and second-order terms above. Since the $M/M/1$ system lives on the non-negative integers, we   add a reflection term at zero to our diffusion. To this end, let us define the  reflected Brownian motion (RBM) satisfying  
 \begin{align}
 Y(t) = Y(0) + \delta(\lambda - \mu) t + \delta \sqrt{\lambda + \mu} W(t) + R(t), \label{eq:rbm}
 \end{align}
 where $R(t)$ is the unique, continuous, and non-decreasing process such that $Y(t) \geq 0$, $R(0) = 0$  and $R(t)$ increases only at those times $t$ when $Y(t) = 0$. Let $Y$ be a random variable having the stationary distribution of this RBM. It is well known that $Y$ is exponentially distributed with mean  $\frac{\delta (\lambda + \mu)}{2( \mu-\lambda) }$,  and so using Lemma 5.2 of \cite{Ross2011} with $Z = Y/\E Y$ there implies 
 \begin{align}
 \E\big( \delta(\lambda - \mu) f'(Y) + \frac{1}{2} \delta^{2} (\lambda + \mu) f''(Y)  \big) + f'(0)  \delta(\mu - \lambda) = 0 \label{eq:rbmbar}
 \end{align} 
 for all $f \in C^{2}(\R_+)$ with $\E \abs{f''(Y)} < \infty$.  One can also derive \eqref{eq:rbmbar} using Ito's lemma for RBMs from Theorem 2 in \cite{HarrReim1981}. 
% Namely, for any $f \in C^{2}(\R_+)$,  
% \begin{align*}
% \E_{x} f(Y(t)) - \E_{x} f(Y(0)) = \E_{x} \Big[ \int_{0}^{t} \big(\delta(\lambda - \mu) f'(Y(s)) + \frac{1}{2} \delta^{2} (\lambda + \mu) f''(Y(s)) \big)  ds + f'(0) R(t) \Big].
% \end{align*}
% Let $Y$ be a random variable having the stationary distribution of this RBM. It is well known that $Y$ is exponentially distributed with mean  $\frac{\delta (\lambda + \mu)}{2( \mu-\lambda) }$.  Picking $f(x) = x$ and initializing $Y(0) \stackrel{d}{=} Y$ above  yields 
%\begin{align*}
% \E \big( R(1) | Y(0) \stackrel{d}{=} Y \big) =  \delta(\mu - \lambda),
%\end{align*} 
% so for any function $f$ such that $\E \abs{f(Y)} < \infty$ and  
%\begin{align*}
%\E \big| \delta(\lambda - \mu) f'(Y) + \frac{1}{2} \delta^{2} (\lambda + \mu) f''(Y) \big| < \infty,
%\end{align*} 
%we can invoke the Fubini-Tonelli theorem to conclude that 

Assume we know \eqref{eq:rbmbar} is satisfied when $f(x) = A f_{h}(x)$, a fact that will be verified by Proposition~\ref{lem:mm1dif} of the following section. Taking expected values   with respect to $Y$ in  \eqref{eq:rbmgx}, and using the fact that $A G_{X} f_h(x) = \E h(X) - A h(x)$,  we   arrive at 
 \begin{align}
& \E h(X) - \E A h(Y) \notag \\
=&\ \E A G_{X} f_h(Y) - \Big( \E\big( \delta(\lambda - \mu) (A f_{h})'(Y) + \frac{1}{2} \delta^{2} (\lambda + \mu) (A f_{h})''(Y)  \big) + (A f_{h})'(0)  \delta(\mu - \lambda) \Big) \notag \\
  =&\ \frac{1}{6} \delta^{3} \E \big( \lambda (A f_h)'''(\xi_{1}(Y)) - \mu (A  \widehat f_h)'''(\xi_{-1}(Y))\big)    - (A f_{h})'(0)  \delta(\mu - \lambda). \label{eq:errorrbm}
\end{align} 

Note that the only term that depends on how we chose our extension $\widehat f_h(\delta k)$  is $(A  \widehat f_h)'''(\xi_{-1}(Y))$ when $\xi_{-1}(Y) \in [-\delta,0]$. Theorem~\ref{thm:1dinterpolant_def} tells us that the upper bound on this term depends on $\Delta^3 \widehat f_h(-\delta)$ and $\Delta^3 \widehat f_h(0) = \Delta^3 f_h(0)$.  A straightforward calculation shows that    $\Delta^3 \widehat f_h(-\delta) = \Delta^2 f_h(0) - \Delta f_h(0)$ because we chose  $\widehat f_h(-\delta) = f_h(0)$. Although $\Delta^2 f_h(0) - \Delta f_h(0)$ looks like a second-order Stein factor, we will see in the next section that it is of the same order of magnitude as $\Delta^3 f_h(0)$, meaning it behaves likes a third-order Stein factor. This is important because as we will see in Proposition~\ref{lem:mm1dif}, third-order Stein factors are smaller than second-order factors by a factor of $\delta$. We also note that a different choice for $\widehat f_h(-\delta)$ would most likely make $\Delta^3 \widehat f_h(-\delta)$ larger than a third-order Stein factor.

\section{Stein Factor Bounds for the $M/M/1$ System} \label{sec:diffbounds}
In this section we show how to use synchronous couplings to bound the Stein factors of the  $M/M/1$ system. We initialize several copies of the same CTMC, each of which have slightly perturbed initial conditions, and we observe how the CTMCs evolve jointly until the time they couple.  As we will see in Section~\ref{sec:synchronous}, the magnitude of the Stein factors depends on a) the coupling time and b) on the distance of the coupled CTMCs relative to  each other before coupling. As such, Stein factors simply measure the sensitivity of the CTMC to perturbations of its initial condition. As to the generality of this approach, one expects it to work well when there is insight into the joint evolution and coupling time of the perturbed and non-perturbed chains.  
The main result of this section is Proposition~\ref{lem:mm1dif} below, which states the relevant Stein factor bounds. We use this lemma  to close the loop on our $M/M/1$ example by proving Theorem~\ref{thm:mm1}.  After discussing synchronous couplings, we also show how the Poisson equation can be used to establish tightness in Section~\ref{sec:tight}.

Our entire discussion relies on our yet-unproven claim that 
\begin{align*}
f_h(\delta k) = \int_{0}^{\infty} \big( \E_{X(0) = \delta k} h(X(t)) - \E h(X) \big) dt 
\end{align*}  
  solves the Poisson equation,  so we now verify this fact.   
\begin{lemma} \label{lem:genpoisson}
Consider a CTMC taking values on a set $E \subset \delta \Z^{d}$ with generator $G_{X}$ given in \eqref{eq:generalctmc}. For $\delta k \in E$, let  $g(\delta k) = \int_{0}^{\infty} \big( \E_{\delta k} h(X(t)) - \E h(X) \big) dt$ and assume that 
\begin{align} 
  g(\delta k)    \text{ is finite for all } \delta k \in E, \label{eq:barbexponential}
\end{align}
 and that $\sum_{\ell \in \Z^{d}} \abs{\beta_{\ell}(\delta k) (g(\delta (k+\ell))-g(\delta k))} < \infty$ for $ \delta k \in E$. 
Then $G_X g(\delta k) = \E h(X) - h(\delta k)$ for all $\delta k \in E$. 
\end{lemma}
Lemma~\ref{lem:genpoisson} is proved in Section~\ref{sec:proofgenpoisson} using an argument similar to one used in  \cite{Barb1988}. In practice, there are several ways to verify that \eqref{eq:barbexponential} holds. One way is by showing that   $\{X(t)\}$ is $h$-exponentially ergodic; i.e.,\ $\abs{\E_{\delta k} h(X(t)) - \E h(X)} \leq c_1 e^{-c_2 t} $ for some $c_1,c_2 > 0$. This is automatically true when the state space $E$ is finite but when $E$ is infinite, the usual way to prove this would be to find a Lyapunov function $V(\delta k)$ such that $G_{X} V(\delta k) \leq -c V(\delta k) + \bar c 1(k \in K)$ for some compact set $K$ and some constants $c,\bar c > 0$. We refer the reader to \cite{MeynTwee1993b} for more on exponential ergodicity. 
 
We now discuss another way  to verify \eqref{eq:barbexponential}   using synchronous couplings.  Note that 
 \begin{align*}
\Big| \int_{0}^{\infty} \big( \E_{\delta k} h(X(t)) - \E h(X) \big) dt\Big|  \leq&\ \int_{0}^{\infty}\sum_{j\in \Z^{d}}\Prob(X = \delta j)  \big|   \E_{\delta k} h(X(t)) - \E_{\delta j} h(X(t))  \big| dt\\
 =&\ \sum_{j\in \Z^{d} }\Prob(X = \delta j) \int_{0}^{\infty} \big|   \E_{\delta k} h(X(t)) - \E_{\delta j} h(X(t))  \big| dt,
 \end{align*}
 where the last equality follows from by Fubini-Tonelli.
%If we find a function $U(\delta k, \delta j)$ such that 
%\begin{align}
%\int_{0}^{\infty} \big|   \E_{\delta k} h(X(t)) - \E_{\delta j} h(X(t))  \big| \leq U(\delta k, \delta j), \quad \text{ and } \quad \E U(\delta k, X) < \infty, \quad \text{ for all } k \in \Z^{d}, \label{eq:difcondition}
%\end{align}
%it will imply \eqref{eq:barbexponential}. As we show below, finding such a  $U(\cdot)$ is naturally tied to using synchronous couplings to establish finite difference bounds. This approach is much closer to the spirit of this paper thanverifying exponential ergodicity.  
 Let us use synchronous couplings to show the right-hand side is finite for the $M/M/1$ model.

Recall the $M/M/1$ generator $G_{X}$ introduced in \eqref{eq:mm1gen} of Section~\ref{sec:bounded}. For the remainder of the section, we let $\{X(t)\}$ and  $X$ represent the corresponding CTMC and stationary distribution, respectively.   Similarly, we let $Y$ have the stationary distribution of the RBM given by  \eqref{eq:rbm} of the same section. Suppose we have proved that for $h \in \dlipone$, 
\begin{align}
\int_{0}^{\infty} \big|   \E_{\delta (k+1)} h(X(t)) - \E_{\delta k} h(X(t))  \big| \leq  \frac{\delta(k+1)}{\mu - \lambda}, \quad k \in \N. \label{eq:diftover}
\end{align}
Using a telescoping sum and the triangle inequality, 
\begin{align*}
& \sum_{j\in \Z}\Prob(X = \delta j) \int_{0}^{\infty} \big|   \E_{\delta k} h(X(t)) - \E_{\delta j} h(X(t))  \big| dt \\
\leq&\  \sum_{j\in \Z}\Prob(X = \delta j) \sum_{i= k \wedge j}^{ k \vee j-1} \int_{0}^{\infty} \big|   \E_{\delta (i+1)} h(X(t)) - \E_{\delta i} h(X(t))  \big| dt.
\end{align*}
Applying \eqref{eq:diftover}, we bound this by 
\begin{align*}
 \sum_{j\in \Z}\Prob(X = \delta j)  \sum_{i= k \wedge j}^{ k \vee j-1}   \frac{\delta(i+1)}{\mu - \lambda} \leq   \sum_{j\in \Z}\Prob(X = \delta j) \frac{\delta(k + j +1)}{\mu - \lambda}(k + j).
\end{align*}
The right-hand side is finite because $\E X^{2} < \infty$, meaning \eqref{eq:barbexponential} is satisfied. The following result confirms \eqref{eq:barbexponential} and presents several Stein factor bounds. It is proved in Section~\ref{sec:synchronous}. 
\begin{proposition} \label{lem:mm1dif}
For any $h \in \dlipone$, \eqref{eq:diftover} holds. Furthermore, for each $ 1 \leq a \leq 3$, if $\abs{\Delta^{v} h(\delta k)} \leq \delta^v$ for all $1 \leq v \leq a$, then 
\begin{align}
 \abs{\Delta^a f_h(\delta k)}  \leq  \frac{\delta^{a}(k+1)}{\mu - \lambda } + \delta^{a-1} \frac{1}{\mu - \lambda}, \quad k \in \N. \label{eq:d3}   
\end{align}  
Lastly, for any function $h : \delta \N \to \R$, 
\begin{align*}
\big( \Delta^{2} f_h(0) - \Delta f_h(0) \big) = \frac{\lambda + \mu }{\mu } \Delta^{3} f_h(0) - \frac{1}{\mu } \Delta^{3} h(0)   - \frac{\lambda}{ \mu } \Delta^{3} f_h(\delta).
\end{align*}
\end{proposition}
The last claim above says that $( \Delta^{2} f_h(0) - \Delta f_h(0) )$ behaves like a third-order Stein factor. Let us  compare Proposition~\ref{lem:mm1dif} to existing Stein factor bounds for the geometric distribution.

 As mentioned in the introduction, the first paper to bound Stein factors for the geometric distribution is \cite{peko1996}. That paper works with indicator test functions  of the form $h(\delta k) = 1(k \in B)$ for $B \subset \N$ instead of allowing $h \in \dlipone$. Bounds for more general test functions can be found in Theorem 1.4 of \cite{Daly2008}. Setting  $q$ and $h(k)$ there to equal $\rho = \lambda/\mu$ and $h(\delta k)/\mu$, respectively, we get the bound 
\begin{align}
  \abs{\Delta^3 f_h(\delta k)}  \leq  \frac{2\delta}{\lambda}, \quad k \in \N, h \in \dlipone. \label{eq:d3alt}   
\end{align} 
Compared to   Proposition~\ref{lem:mm1dif}, the right-hand side of \eqref{eq:d3alt} is  independent of   $k$, and only requires $h \in \dlipone$, making it more convenient to work with.   However, unlike synchronous couplings, the proof of \eqref{eq:d3alt} in \cite{Daly2008} is much harder to generalize to multidimensional settings.  We now state and prove a bound on the approximation error between $Y$ and $X$. 
\begin{theorem}
\label{thm:mm1}
 There exists a constant $C > 0$ such that for $\rho < 1$ and $h^{*} \in \lipone$, 
\begin{align*}
\abs{\E h^{*}(X) - \E h^{*}(Y) } \leq C \delta \Big(1 + \frac{1}{\rho} \Big).
\end{align*}
\end{theorem}
To prove Theorem~\ref{thm:mm1} we will use  the third-order bound in \eqref{eq:d3alt}  instead of   Proposition~\ref{lem:mm1dif}.  After the proof,  we comment on how to get a comparable bound using Proposition~\ref{lem:mm1dif}. 
Before proving the theorem, let us say a few words on the possible choices of $\delta$. It is well known that $X \approx Y$ when $\rho \approx 1$, and that   $\E X = \delta \rho/(1-\rho)$. Choosing $\delta = 1$, Theorem~\ref{thm:mm1} tells us that the approximation error $\abs{\E X - \E Y}$ does not grow even though $\E X \to \infty$ as $\rho \to 1$. However, since both $X$ and $Y$ diverge, we cannot conclude that $X$ converges to $Y$. To ensure convergence, we recall that $\E Y = \frac{\delta (\lambda + \mu)}{2( \mu-\lambda) }$. Choosing $\delta = (1-\rho)$  ensures that $\{X\}_{\rho < 1}$ and $\{Y\}_{\rho < 1}$ are  tight, and Theorem~\ref{thm:mm1} then implies that $X$ converges to $Y$ in distribution as $\rho \to 1$. As discussed in the introduction, tightness of the prelimit sequence is a sought-after property  because, when combined with process-level convergence to some diffusion limit, tightness implies convergence of stationary distributions as well.  We discuss in Section~\ref{sec:tight} below how one can use the Poisson equation to establish tightness.  

\startproof{Proof of Theorem~\ref{thm:mm1}}
Let $h: \delta \N \to \R$ be the restriction of $h^{*}(x)$ to $\delta \N$.  
Since $Y$ is exponentially distributed,  Proposition~\ref{lem:mm1dif} implies \eqref{eq:rbmbar} holds with $f(x) = A  f_h(x)$ there. 
%\begin{align*}
% \E\big( \delta(\lambda - \mu) ( A  f_h)'(Y) + \frac{1}{2} \delta^{2} (\lambda + \mu) ( A  f_h)''(Y)  \big) + ( A  f_h)'(0)  \delta(\mu - \lambda) = 0.
%\end{align*}
 Consequently,   \eqref{eq:errorrbm} holds, which we recall below:
\begin{align*}
 \E h(X) - \E Ah(Y) =&\ \frac{1}{6} \delta^{3} \E \big( \lambda (A f_h)'''(\xi_{1}(Y)) - \mu (A  \widehat f_h)'''(\xi_{-1}(Y))\big)   - (A f_{h})'(0)  \delta(\mu - \lambda).
\end{align*}
Once we bound the right-hand side above, Lemma~\ref{lem:convdet} will imply the theorem. Inequality \eqref{eq:multibound} from Theorem~\ref{thm:1dinterpolant_def} and the Stein factor bound in \eqref{eq:d3alt} imply
\begin{align*}
 \delta^{3} \lambda \abs{(A f_h)'''(\xi_{1}(Y))} \leq&\ C \lambda  \max_{0\leq i \leq 1} \abs{\Delta^{3} f_h\big(\delta (k(\xi_1(Y))  + i)\big) } \leq C  \delta.
\end{align*}
Similarly, the  bounds from Proposition~\ref{lem:mm1dif} imply $\abs{(A f_{h})'(0)  \delta(\mu - \lambda)} \leq  C\delta^2$. Lastly
\begin{align*}
\delta^{3}  \mu \abs{(A  \widehat f_h)'''(\xi_{-1}(Y))} \leq&\ C \mu  \max_{0\leq i \leq 1} \abs{\Delta^{3}  \widehat f_h\big(\delta (k(\xi_{-1}(Y))  + i)\big) }.
\end{align*}
As discussed at the end of Section~\ref{sec:bounded}, $\Delta^{3}  \widehat f_h\big(\delta (k(\xi_{-1}(Y))  + i) = \Delta^{3}   f_h\big(\delta (k(\xi_{-1}(Y))  + i)$ if $\xi_{-1}(Y) + i  \geq 0$, and otherwise it equals $\Delta^{2} f_h(0) - \Delta f_h(0)$. Using the form of $\Delta^{2} f_h(0) - \Delta f_h(0)$ from Proposition~\ref{lem:mm1dif} together with \eqref{eq:d3alt}, we get 
\begin{align*}
C \mu  \max_{0\leq i \leq 1} \abs{\Delta^{3}  \widehat f_h\big(\delta (k(\xi_{-1}(Y))  + i)\big) } \leq C \delta \frac{1}{\rho}.
\end{align*}
\finishproof
The  approximation error can also be bounded using only the bounds from   Proposition~\ref{lem:mm1dif} instead of \eqref{eq:d3alt}. Consider the error term 
\begin{align*}
 \delta^{3} \lambda \E \abs{(A f_h)'''(\xi_{1}(Y))} \leq&\ C \lambda  \E \max_{0\leq i \leq 1} \abs{\Delta^{3} f_h\big(\delta (k(\xi_1(Y))  + i)\big) } 
\end{align*}
from the proof of the theorem above. If we apply Proposition~\ref{lem:mm1dif} and the fact that $\delta k(\xi_1(Y)) \leq Y+\delta $, we get 
\begin{align*}
C \lambda  \E \max_{0\leq i \leq 1} \abs{\Delta^{3} f_h\big(\delta (k(\xi_1(Y))  + i)\big) }  \leq C \lambda \delta^{2}\frac{\E Y + 1}{\mu - \lambda} =&\ C \lambda \delta^{3}\frac{\lambda + \mu}{2(\mu - \lambda)^2} + C \lambda \delta^{2}\frac{1}{\mu - \lambda} \\
=&\ C \rho \delta^{3}\frac{\rho + 1}{2(1-\rho)^2} + C \rho  \delta^{2}\frac{1}{1-\rho}.
\end{align*}
In the second-last equality we used $\E Y = \frac{\delta (\lambda + \mu)}{2( \mu-\lambda) }$. Choosing $\delta = (1-\rho)$ means the term on the right-hand side is bounded by $C \delta$, giving a comparable bound to the one in Theorem~\ref{thm:mm1}. Other error terms can be bounded similarly.  As discussed previously, the choice of $\delta = (1-\rho)$ is natural because it ensures tightness of $\{X\}_{\rho < 1}$ and $\{Y\}_{\rho < 1}$.
%\begin{remark}
%\label{rem:errs}
%We return to our postponed discussion of $\E \varepsilon_h(Y)$ and $\E \varepsilon_f(Y)$ from Proposition~\ref{lem:interchangebounded}. 
%Suppose $\delta = (1-\rho)$.  Let us discuss the error terms $\E \varepsilon_h(Y)$ and $\E \varepsilon_f(Y)$ and 
%The magnitude of the error term was dictated by the error terms involving $(A    f_h)'''(\xi_{1}(Y))$ and $(A  \widehat f_h)'''(\xi_{-1}(Y))$. The error terms $\E \varepsilon_h(Y)$ and $\E \varepsilon_f(Y)$ could have been 
%\end{remark} 

\subsection{Synchronous Couplings}
\label{sec:synchronous}
We now use synchronous couplings to prove Proposition~\ref{lem:mm1dif}. 

\startproof{Proof of Proposition~\ref{lem:mm1dif}}

\textbf{First-order factors.}
Consider two $M/M/1$ systems, whose customer counts (scaled by $\delta$) are $\{X^{(0)}(t)\}$ and $\{X^{(1)}(t)\}$. We refer to these as system $0$ and system $1$, respectively. We couple the two systems by setting $X^{(1)}(0) = X^{(0)}(0)+1$  and defining their joint evolution via the following transition rate table.
\begin{table}[h!]
\caption{Transitions of the joint process $\{(X^{(1)}(t), X^{(0)}(t))\}$ in  state $   (x^{(1)},x^{(0)}) $. }
\centering
\label{tab:mm1}
\begin{tabular}{|c|c|c|}
\hline
\# &Rate & Transition      \\ \hline
1& $ \lambda  $ & $(x^{(1)}+ \delta,x^{(0)}+ \delta)$  \\ \hline
2&$\mu 1(x^{(0)}>0)  $ & $(x^{(1)}- \delta,x^{(0)}- \delta)$  \\ \hline
3&$\mu 1(x^{(1)}>0,x^{(0)}=0)  $ & $(x^{(1)}- \delta,x^{(0)})$  \\ \hline
\end{tabular}
\end{table}

The transition table can be interpreted as follows. Both systems have the same customer arrival stream. The customers present in system 0 at time $t = 0$, and all newly arriving customers, have an identical counterpart in system 1. The only difference between the two is the extra initial customer in system 1, who behaves like a low-priority customer that only gets served when there are no other customers, and is preempted by new arrivals. The two systems couple  once this extra customer is served. We refer to systems 1 and 0 as a synchronous coupling  because the two systems are driven by the same underlying stochastic processes, i.e.,\ arrivals and services.

To bound $\Delta f_h(\delta k)$, define $\tau^{(i)}(\delta k) = \inf_{t \geq 0} \{X^{(i)}(t) = \delta k\}$. The discussion above implies
\begin{align}
\Delta f_h(\delta k)  =&\ \int_{0}^{\infty} \big(\E_{\delta (k+1)}  h(X (t)) -  \E_{\delta k}h(X (t)) \big) dt  \notag \\
=&\ \int_{0}^{\infty} \E_{X^{(0)}(0)=\delta k}   \Big( h(X^{(1)}(t)) -  h(X^{(0)}(t)) \Big)  dt \notag \\
=&\ \int_{0}^{\infty} \E_{X^{(0)}(0)=\delta k} \bigg[ 1 (t \leq  \tau^{(1)}(0)  )   \Big( h(X^{(0)}(t)+\delta) -  h(X^{(0)}(t)) \Big)  \bigg] dt.  \label{eq:dfk_mm1}
\end{align}
We emphasize that the last equality above is true because systems 0 and 1 always maintain a constant gap of a single customer until they couple. We bound $\E_{X^{(0)}(0)=\delta k}  \tau^{(1)}(0)$ by combining the Lyapunov function $V(\delta k ) = k$ with Dynkin's formula. Observe that $V(\delta k)$ satisfies $G_X V(\delta k) =  \lambda - \mu < 0$ whenever $k > 0$, which means that
%\begin{align*}
%\E_{X^{(1)}(0)=\delta (k+1)} V\big(X^{(1)}( \tau^{(1)}(0))  \big) - V(\delta k ) =&\ \E_{X^{(1)}(0)=\delta (k+1)} \int_{0}^{ \tau^{(1)}(0)} (\lambda - \mu) ds,
%\end{align*}
%which implies
\begin{align}
 \E_{X^{(0)}(0)=\delta k}  \tau^{(1)}(0)  =   \frac{\E_{X^{(1)}(0)=\delta (k+1)} \big(V(X^{(1)}(0) ) -   V(X^{(1)}( \tau^{(1)}(0) ) )\big)}{\mu - \lambda} = \frac{k+1}{\mu - \lambda}. \label{eq:taubound}
\end{align}  
To justify the above equality, we refer the reader to  the proof of Theorem 4.3.i of \cite{MeynTwee1993b}, which is a direct application of Dynkin's formula.
%To be precise, the inequality above is actually  an equality because $ V(X^{(1)}( \tau^{(1)}(0)) ) = 0$, but this does not matter to us.  
Combining \eqref{eq:dfk_mm1},  \eqref{eq:taubound} and the fact that $h \in \dlipone$ proves
\begin{align*}
\abs{\Delta f_h(\delta k)} \leq \frac{\delta (k+1)}{\mu - \lambda}, \quad k \in \N.
\end{align*}
In fact, we have proved the stronger statement  \eqref{eq:diftover}.

\begin{figure}[h!]
\centering
\includegraphics[scale=1]{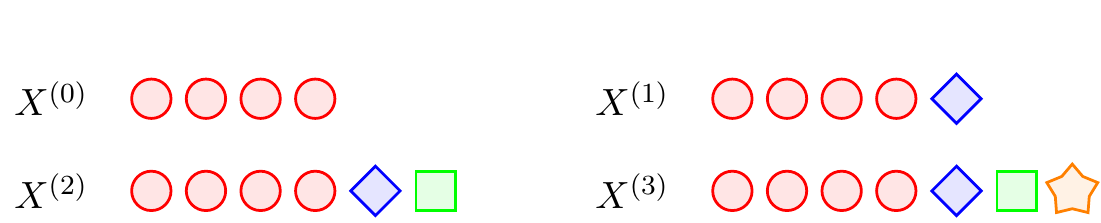}
\caption{The initial state of systems 0,1,2,3 when system 0 starts with 4 customers.  The diamonds, squares, and stars  represent the extra customers. }
\label{fig:starting_mm1_4}
\end{figure}

\textbf{Second- and third-order factors.} 
We now prove the third-order bounds. Second-order bounds follow analogously.  In addition to systems 0 and 1, we let $\{X^{(2)}(t)\}$ and $\{X^{(3)}(t)\}$ represent systems 2 and 3. System 2 is an identical copy of system 1 with one additional low-priority customer, and system 3 is a copy of system 2 with yet another low-priority customer. The relationship between the four systems is visualized in Figure~\ref{fig:starting_mm1_4}, where we note that $X^{(3)}(0) =X^{(2)}(0)+\delta = X^{(1)}(0)+2\delta = X^{(0)}(0)+3\delta$.  The transitions of the joint chain are  formally  defined in Table~\ref{tab:mm1d2} below.
 
\begin{table}[h!]
\caption{Transitions of the joint process $\{(X^{(3)}(t),X^{(2)}(t),X^{(1)}(t), X^{(0)}(t))\}$ in  state $   (x^{(3)},x^{(2)},x^{(1)},x^{(0)}) $. }
\centering
\label{tab:mm1d2}
\begin{tabular}{|c|c|c|}
\hline
\# &Rate & Transition      \\ \hline
1& $ \lambda  $ & $(x^{(3)}+ \delta,x^{(2)}+ \delta,x^{(1)}+ \delta,x^{(0)}+ \delta)$  \\ \hline
2&$\mu 1(x^{(0)}>0)  $ & $(x^{(3)}- \delta,x^{(2)}- \delta,x^{(1)}- \delta,x^{(0)}- \delta)$  \\ \hline
3&$\mu 1(x^{(1)}>0,x^{(0)}=0)  $ & $(x^{(3)}- \delta,x^{(2)}- \delta,x^{(1)}- \delta,x^{(0)})$  \\ \hline
4&$\mu 1(x^{(2)}>0,x^{(1)}=0)  $ & $(x^{(3)}- \delta,x^{(2)}- \delta,x^{(1)},x^{(0)})$  \\ \hline
5&$\mu 1(x^{(3)}>0,x^{(2)}=0)  $ & $(x^{(3)}- \delta,x^{(2)},x^{(1)},x^{(0)})$  \\ \hline
\end{tabular}
\end{table}
%Observe that systems 0 and 1 are identical for all $t \geq  \tau^{(1)}(0)$.

\noindent It follows that
%\begin{align}
%\Delta^3 f_h(\delta k)  =&\ \int_{0}^{\infty} \E_{X^{(0)}(0)=\delta k}   \Big( h(X^{(3)}(t)) - 3 h(X^{(2)}(t)) + 3 h(X^{(1)}(t)) -  h(X^{(0)}(t)) \Big)  dt \notag \\
%=&\ \int_{0}^{\infty} \E_{X^{(0)}(0)=\delta k} \bigg[ 1 (t \leq  \tau^{(1)}(0)  ) \Delta^3 h(X^{(0)}(t)) \bigg] dt \notag \\
%& + \int_{0}^{\infty} \E_{X^{(1)}(0)=0}   \Big(  h(X^{(3)}(t)) - 3 h(X^{(2)}(t)) + 2 h(X^{(1)}(t)) \Big)  dt  \notag  \\
%=&\ \int_{0}^{\infty} \E_{X^{(0)}(0)=\delta k} \bigg[ 1 (t \leq  \tau^{(1)}(0)  ) \Delta^3 h(X^{(0)}(t)) \bigg] dt + \Delta^{2} f_h(0) - \Delta f_h(0) \notag \\
%=&\ \int_{0}^{\infty} \E_{X^{(0)}(0)=\delta k} \bigg[ 1 (t \leq  \tau^{(1)}(0)  ) \Delta^3 h(X^{(0)}(t)) \bigg] dt \notag \\
%&+ \int_{0}^{\infty} \E_{X^{(0)}(0)= 0} \bigg[ 1 (t \leq  \tau^{(1)}(0)  ) \Delta^2 h(X^{(0)}(t)) \bigg] dt. \label{eq:d3calcs}
%\end{align} 
\begin{align}
\Delta^3 f_h(\delta k)  =&\ \int_{0}^{\infty} \E_{X^{(0)}(0)=\delta k}   \Big( h(X^{(3)}(t)) - 3 h(X^{(2)}(t)) + 3 h(X^{(1)}(t)) -  h(X^{(0)}(t)) \Big)  dt \notag \\
=&\ \int_{0}^{\infty} \E_{X^{(0)}(0)=\delta k} \bigg[ 1 (t \leq  \tau^{(1)}(0)  ) \Delta^3 h(X^{(0)}(t)) \bigg] dt \notag \\
& + \int_{0}^{\infty} \E_{X^{(1)}(0)=0}   \Big(  h(X^{(3)}(t)) - 3 h(X^{(2)}(t)) + 2 h(X^{(1)}(t)) \Big)  dt  \notag  \\
=&\ \int_{0}^{\infty} \E_{X^{(0)}(0)=\delta k} \bigg[ 1 (t \leq  \tau^{(1)}(0)  ) \Delta^3 h(X^{(0)}(t)) \bigg] dt + \Delta^{2} f_h(0) - \Delta f_h(0).   \label{eq:d3calcs}
\end{align} 
Similarly,
\begin{align*}
 \Delta^2 f_h(\delta k)  =&\ \int_{0}^{\infty} \E_{X^{(0)}(0)=\delta k}   \Big( h(X^{(2)}(t)) - 2 h(X^{(1)}(t)) -  h(X^{(0)}(t)) \Big)  dt \notag \\
=&\ \int_{0}^{\infty} \E_{X^{(0)}(0)=\delta k} \bigg[ 1 (t \leq  \tau^{(1)}(0)  ) \Delta^2 h(X^{(0)}(t)) \bigg] dt + \Delta f_h(0).
\end{align*}
Combining the representations for $\Delta^3 f_h(\delta k) $ and $\Delta^2 f_h(\delta k) $  above,  we conclude that 
\begin{align*}
\Delta^3 f_h(\delta k)   =&\ \int_{0}^{\infty} \E_{X^{(0)}(0)=\delta k} \bigg[ 1 (t \leq  \tau^{(1)}(0)  ) \Delta^3 h(X^{(0)}(t)) \bigg] dt \notag \\
&+ \int_{0}^{\infty} \E_{X^{(0)}(0)= 0} \bigg[ 1 (t \leq  \tau^{(1)}(0)  ) \Delta^2 h(X^{(0)}(t)) \bigg] dt.  
\end{align*} 
Assuming $h \in \dlipone$,  $\abs{\Delta^{2} h(\delta k)} \leq \delta^2$, and $\abs{\Delta^{3} h(\delta k)} \leq \delta^3$, we apply \eqref{eq:taubound} to conclude that
\begin{align*}
\abs{\Delta^{3} f_h(\delta k)} \leq  \frac{\delta^{3}(k+1)}{\mu - \lambda } + \delta^{2} \frac{1}{\mu - \lambda}.
\end{align*} 
The bound on $\abs{\Delta^{2} f_h(\delta k)}$ is obtained similarly. Lastly, if we consider the first equality in \eqref{eq:d3calcs} with $k = 0$, and consider what happens after the first jump of the CTMC, we get
\begin{align*}
 \Delta^{3} f_h(0) = \frac{1}{\lambda + \mu} \Delta^{3} h(0) + \frac{\mu }{\lambda + \mu }\big( \Delta^{2} f_h(0) - \Delta f_h(0) \big) + \frac{\lambda}{\lambda + \mu } \Delta^{3} f_h(\delta).
\end{align*}
\finishproof

\subsection{Tightness via the  Poisson Equation}
\label{sec:tight}
One way to show a sequence of CTMC stationary distributions is tight is by bounding  $\E \abs{X}$. For the $M/M/1$ system this task is trivial because the stationary distribution is known. However, obtaining a useful upper bound on $\E \abs{X}$ is   harder for more complicated systems and usually involves using some kind of Lyapunov function. We can also use the Poisson equation 
\begin{align*}
G_X f_h(\delta k) = \E h(X) - h(\delta k), \quad k \geq 0 
\end{align*}
as follows.  
Pick $h(\delta k) = \abs{\delta k}$ and evaluating the above at $k = 0$ to get 
\begin{align*}
\lambda \big( f_h(\delta) - f_h(0) \big) = \lambda \Delta f_h(0)  = \E \abs{X} = \E X.
\end{align*}
Recall \eqref{eq:dfk_mm1} and \eqref{eq:taubound} from the previous section,  which  imply  that  $  \E X = \delta \lambda/(\mu -\lambda ) = \delta \rho /(1-\rho)$. Choosing $\delta$ to be any constant multiple of $1-\rho$ ensures  that $\{X\}_{\rho < 1}$ is tight.

The main takeaway is that the problem of tightness is equivalent to bounding $G_X f_h(\delta k)$ at \emph{a single point}, usually the fluid equilibrium of the CTMC. At the fluid equilibrium, $G_{X} f_h(\delta k)$ typically consists of some combination of first- and second-order Stein factors.

\section{Misalignment of Diffusion Synchronous Couplings}
\label{sec:compare_couplings}
We have presented the prelimit approach as a parallel to the diffusion approach for the purposes of bounding $\abs{\E h(X) - \E h(Y)}$. As we have seen, the  main challenge with either approach is bounding the differences/derivatives of the respective Poisson equation solution. In theory, any bound achievable using one approach should be achievable with the other. In practice, there are slight technical differences between working with a discrete-valued CTMC and a diffusion living on a continuum. In this section we illustrate one  technical nuance that arises when we use synchronous couplings to bound the third derivatives of $f_{h^{*}}(x)$ in the diffusion approach. We term this the ``misalignment of synchronous couplings''. The main takeaway is that the misaligned synchronous couplings add extra complexity to the problem. In contrast, the analogous analysis using the prelimit approach in Section~\ref{sec:synchronous} is cleaner because the CTMC is restricted to the grid.

Recall  the generic diffusion process $\{Y(t) \in \R \}$  defined in \eqref{eq:sdediffusion}. We assume for simplicity that the diffusion coefficient $a(x) = a$ for all $x \in \R$  and define the synchronous couplings 
\begin{align*}
Y^{(i)}(t) = Y^{(0)}(0) + i \varepsilon + \int_{0}^{t} b(Y^{(i)}(s)) ds + \int_{0}^{t} \sqrt{a} d W(s), \quad i = 0,1,2,3.
\end{align*}
The four couplings start at different initial conditions  but share the same Brownian motion.  Since  $f_{h^{*}}(x)$ is given by \eqref{eq:diffusionpoissonsolution}, it follows that  
\begin{align}
&\frac{\partial^3}{\partial x^3} f_{h^{*}}(x) \notag \\
 =&\ \lim_{\varepsilon \to 0} \frac{1}{\varepsilon^{3}} \int_{0}^{\infty}  \E_{Y^{(0)}(0) = x  } \Big(h^{*}(Y^{(3)}(t)) - 3  h^{*}(Y^{(2)}(t))  + 3  h^{*}(Y^{(1)}(t))  -   h^{*}(Y^{(0)}(t))   \Big) dt. \label{eq:misalignstart}
\end{align}
To show   the integral on the right-hand side is finite, one must characterize the speed at which the synchronous couplings converge to one another. Furthermore, the integral  must be of order $\varepsilon^{3}$   for the limit to exist. Let us consider this last point further. 

Given a sufficiently differentiable function $g: \R \to \R$,  we know that its derivatives can be approximated by finite differences. For instance, Taylor expansion tells us that
\begin{align}
 g'''(x) \approx \frac{\big( g(x''') - 3 g(x'') + 3 g(x') - g(x)\big) }{\varepsilon^{3}} \label{eq:stencil}
\end{align}
when $x''' = x+ 3\varepsilon$, $x'' = x + 2 \varepsilon$, and $x' = x + \varepsilon$. The precise spacing of $x,x',x'',x'''$ relative to each other is essential for the limit (as $\varepsilon \to 0$) of the right-hand side in \eqref{eq:stencil} to exist. For example, if $x''' = x + 4\varepsilon$, then the numerator is now of order $\varepsilon$ instead of $\varepsilon^3$, and the right-hand side diverges  as $\varepsilon \to 0$.  Therefore, one way to show that the integral in \eqref{eq:misalignstart} is of order $\varepsilon^{3}$ is to prove that the diffusion couplings maintain the appropriate spacing relative to each other so that the integrand is of order $\varepsilon^3$ for each $t \geq 0$. 

Indeed, this is precisely the result of Lemma 3.3 of \cite{GorhMack2016}, which says that if $h(x)$ and $b(x)$ are smooth enough,  and if $b(x)$ is  $k$-strongly concave, then  
\begin{align}
\abs{h^{*}(Y^{(3)}(t)) - 3  h^{*}(Y^{(2)}(t))  + 3  h^{*}(Y^{(1)}(t))  -   h^{*}(Y^{(0)}(t)) }  \leq \varepsilon^3 C e^{-kt/2} \label{eq:gormack}
\end{align}
almost surely, where the constant $C$ depends on $k$, $h(x)$ and $b(x)$. The above inequality then implies that 
\begin{align}
 \lim_{\varepsilon \to 0} \frac{1}{\varepsilon^{3}} \int_{0}^{\infty}    \E_{Y^{(0)}(0) = x  } \Big|h^{*}(Y^{(3)}(t)) - 3  h^{*}(Y^{(2)}(t))  + 3  h^{*}(Y^{(1)}(t))  -   h^{*}(Y^{(0)}(t))   \Big| dt \leq 2C/k. \label{eq:gormakimply}
\end{align}
Similarly, \eqref{eq:gormack} also holds for $d$-dimensional diffusions with constant diffusion coefficients. Unfortunately, if the   assumptions on the drift  are violated, e.g.\ the drift is only Lipschitz-continuous or the diffusion has a reflecting boundary, then \eqref{eq:gormack} no longer holds because the diffusion couplings become misaligned. This misalignment complicates the problem  of bounding \eqref{eq:misalignstart} because one cannot use \eqref{eq:gormack} anymore. 

As an example, we now illustrate  how this misalignment occurs in the RBM that approximates the $M/M/1$ system. In contrast, the discrete nature of the prelimit approach prevents this kind of misalignment from happening. Recall the  RBM  defined in \eqref{eq:rbm},
\begin{align*}
Y(t) = Y(0) + \delta (\lambda - \mu)t+  \delta\sqrt{(\lambda + \mu)}W(t) + R(t), \quad t \geq 0,
\end{align*}
and let $Y$ be the random variable having its stationary distribution. It was shown in \cite{HarrReim1981} that 
\begin{align*}
R(t) = -\inf_{0 \leq s \leq t}  \Big\{ Y(0) + \delta (\lambda - \mu)s+  \delta\sqrt{(\lambda + \mu)}W(s) \Big\}. 
\end{align*} 
We wish to bound the third derivative of $f_{h^{*}}(x)$. 
% \begin{align*}
%f_{h^{*}}(x) = \int_{0}^{\infty} \big( \E_{Y(0) = x}  h^{*}(Y(t))  - \E h^{*}(Y) \big) dt, \quad x \geq 0.
%\end{align*}
For simplicity, we choose $h^{*}(x) = x$. Let us define the four  coupled processes 
\begin{align}
Y^{(i)}(t) =&\  Y^{(0)}(0) + i \varepsilon  + \delta (\lambda - \mu)t+  \delta\sqrt{(\lambda + \mu)}W(t) + R^{(i)}(t), \notag \\
 \text{ where } \quad R^{(i)}(t) =&\ -\inf_{0 \leq s \leq t}  \Big\{ Y^{(i)}(0)  + \delta (\lambda - \mu)s+  \delta\sqrt{(\lambda + \mu)}W(s) \Big\}, \quad i = 0,1,2,3. \label{eq:ridef}
\end{align}
%We refer to $\{Y^{(i)}(t)\}$ as diffusion $i$.  
We also define $D_3(t) =Y^{(3)}(t) - 3Y^{(2)}(t) + 3 Y^{(1)}(t) - Y^{(0)}(t)$.  It follows that 
\begin{align}
  \frac{\partial^3}{\partial x^3} f_{h^{*}}(x)  =&\   \lim_{\varepsilon \to 0} \frac{1}{\varepsilon^3}  \int_{0}^{\infty}  \E_{Y^{(0)}(0) = x}D_3(t) dt. \label{eq:d3_starting_point}
\end{align}
%Observe that 
%\begin{align*}
%Y^{(i+1)}(t) - Y^{(i)}(t) =  \varepsilon + R^{(i+1)}(t) - R^{(i)}(t), \quad i = 0,1,2,
%\end{align*}
%and so the integrand in \eqref{eq:d3_starting_point} equals zero as long as the $R^{(i)}(t)$s equal zero, i.e.\ until the first time that diffusion 0 hits the origin. Furthermore, the integrand is zero from the time that the four diffusions couple, which occurs when diffusion 3 hits the origin. 
We define 
\begin{align*}
\gamma_1 = \inf_{t \geq 0} \{Y^{(1)}(t) = 3\varepsilon/4 \}, \quad \gamma_2 =  \inf_{t \geq 0} \{Y^{(1)}(t) =  \varepsilon/4 \}.
\end{align*}
We will  prove at the end of this section that  
\begin{align}
 D_3(t) \leq -\varepsilon/4, \quad \text{ for } t \in [\gamma_1,\gamma_2]. \label{eq:toproverbm}
\end{align}
We see that \eqref{eq:toproverbm} violates \eqref{eq:gormack}. Furthermore, the expected hitting time of a fixed level by a Brownian motion with drift is well known and implies that $\E(\gamma_2 - \gamma_1) = \varepsilon/(2\delta(\mu-\lambda)) $. Therefore, the integral in \eqref{eq:d3_starting_point} equals
\begin{align}
&\frac{1}{\varepsilon^3}  \int_{0}^{\infty}  \E_{Y^{(0)}(0) = x}(D_3(t) 1(t \in[\gamma_1,\gamma_2] )) dt   +\frac{1}{\varepsilon^3} \int_{0}^{\infty}  \E_{Y^{(0)}(0) = x}(D_3(t) 1(t \not\in[\gamma_1,\gamma_2] )) dt , \label{eq:int12}
\end{align}
and the first term is bounded from above by $   -(8\delta(\mu-\lambda) \varepsilon)^{-1}$, which diverges as $\varepsilon \to 0$. Therefore, unlike in \eqref{eq:gormakimply}, we cannot just bound $\abs{D_3(t)}$ and take the limit as $\varepsilon \to 0$.

In reality, $\frac{\partial^3}{\partial x^3} f_{h^{*}}(x)$ is well defined and the right-hand side of \eqref{eq:d3_starting_point} exists and it is possible to prove that the second integral in \eqref{eq:int12} contains a positive term of order $1/\varepsilon$ that  cancels out the first integral. Indeed, Theorem 1.2 in \cite{Daly2008} presents bounds on $\frac{\partial^k}{\partial x^k} f_{h^{*}}(x)$ for any $k \geq 2$.  We conclude by verifying \eqref{eq:toproverbm}. By definition, 
\begin{align*}
Y^{(i+1)}(t) - Y^{(i)}(t) =  \varepsilon + R^{(i+1)}(t) - R^{(i)}(t), \quad i = 0,1,2,
\end{align*}
for  every $t \geq 0$. Since $R^{(i)}(t) = 0$ for $i = 1,2,3$ when $t < \inf_{s \geq 0} \{Y^{(1)}(s) = 0\}$, we have that
\begin{align*}
Y^{(3)}(t) - Y^{(2)}(t) =Y^{(2)}(t) - Y^{(1)}(t) = \varepsilon, \quad t \in [\gamma_1,\gamma_2] 
\end{align*}
because $\gamma_2 < \inf_{s \geq 0} \{Y^{(1)}(s) = 0\}$. Thus, 
\begin{align*}
D_3(t)  =&\ -\varepsilon +  Y^{(1)}(t) - Y^{(0)}(t)=  R^{(1)}(t) - R^{(0)}(t) = - R^{(0)}(t), \quad t \in [\gamma_1,\gamma_2].
\end{align*}
One can check that $R^{(0)}(\gamma_1) = \varepsilon/4$ and $R^{(0)}(\gamma_2) = 3\varepsilon/4$ using the form of $R^{(i)}(t)$ in \eqref{eq:ridef}. Since $R^{(0)}(t)$ is non-decreasing, we have  $R^{(0)}(t) \in [\varepsilon/4, 3\varepsilon/4]$ when $t \in [\gamma_1,\gamma_2]$, proving \eqref{eq:toproverbm}.

 \section{Conclusion}
\label{sec:conclusion}
In this paper we introduced the prelimit generator comparison approach and used the $M/M/1$ model to illustrate it in practice. In applying the approach, we overcame two technical challenges. First, we used an interpolator to extend the prelimit Poisson equation to the continuum. Second, we  interchanged the interpolator with the CTMC generator. Our solution to both challenges extends beyond the $M/M/1$ model. Appendices~\ref{sec:interpolation} and \ref{sec:mdimunbound} contain the general multidimensional  interpolation and interchange results.  These are intended to simplify as much as possible the tedious aspects of the prelimit approach, and to make it easy for readers to apply the approach to their own problem.

One direction we have not considered is working with the Kolmogorov distance. For two real-valued random variables $U,U'$, the Kolmogorov distance is defined as 
\begin{align*}
d_{K}(U,U') = \sup_{z \in \R} \big| \E \big(1(U \geq z)\big) - \E \big(1(U' \geq z)\big)\big|. 
\end{align*}
It is well known (e.g.,\ \cite{BravDaiFeng2016}) that the discontinuity in the test functions $1(\cdot \geq z)$  makes working with the Kolmogorov distance more difficult than the Wasserstein. Even though we deal with discrete functions and their interpolations, the issue with the discontinuity in $1(\cdot \geq z)$ will still come up in the difference bounds on $f_h(x)$. Although there  are undoubtedly technical challenges to overcome, the author believes that  the prelimit approach can be used as a tool to bound the Kolmogorov distance.

% Appendix here
% Options are (1) APPENDIX (with or without general title) or
%             (2) APPENDICES (if it has more than one unrelated sections)
% Outcomment the appropriate case if necessary
%
% \begin{APPENDIX}{<Title of the Appendix>}
% \end{APPENDIX}
%
%   or
%
 \begin{APPENDICES}

\section{Multidimensional Interpolation}
\label{sec:interpolation}
We now prove Theorem~\ref{thm:1dinterpolant_def}. We then state and prove the multidimensional version --  Theorem~\ref{thm:interpolant_def}.
\startproof{Proof of Theorem~\ref{thm:1dinterpolant_def}}

Given $f: K \to \R$,   for each $k \in \Z$ such that $\delta k \in K_{4}$  we define 
\begin{align}
P_{k}(x) =&\ f(\delta k) + \Big(\frac{x-\delta k}{\delta} \Big)(\Delta - \frac{1}{2}\Delta^2 + \frac{1}{3}\Delta^3) f(\delta k)  \notag  \\
&+ \frac{1}{2} \Big(\frac{x-\delta k}{\delta} \Big)^2\big(\Delta^2 - \Delta^3\big) f(\delta k) + \frac{1}{6}  \Big(\frac{x-\delta k}{\delta} \Big)^3 \Delta^3 f(\delta k)\notag \\
& -\frac{23}{3} \Big(\frac{x-\delta k}{\delta} \Big)^4 \Delta^4 f(\delta k)  +\frac{41}{2} \Big(\frac{x-\delta k}{\delta} \Big)^5\Delta^4 f(\delta k)  \notag  \\
&- \frac{55}{3}  \Big(\frac{x-\delta k}{\delta} \Big)^6\Delta^4 f(\delta k)  +\frac{11}{2} \Big(\frac{x-\delta k}{\delta} \Big)^7\Delta^4 f(\delta k) , \quad   x \in \R. \label{eq:pplus}
\end{align}
From \eqref{eq:pplus} we have $P_{k}(\delta k) = f(\delta k)$, implying \eqref{eq:interpolates}. Since   $P_{k}(x) \in C^{\infty}(\R)$, we know $A f(x)$ is infinitely differentiable on $\text{Conv}(K_4) \setminus K_4$. Furthermore, it is straightforward to verify that
\begin{align}
\frac{\partial^{a}}{\partial x^{a}} P_{k-1}(x)\Big|_{x =\delta k}  = \frac{\partial^{a}}{\partial x^{a}} P_{k}(x)\Big|_{x =\delta k} , \quad \text{ for } a = 0,1,2,3. \label{eq:onedsmooth}
\end{align}
The property above implies $A f(x) = P_{k(x)}(x) \in C^{3}(\text{Conv}(K_4))$.  The  weights $\alpha^{k}_{k+i}(x)$  can be read off by combining the coefficients corresponding to $f(\delta (k+i))$ in  \eqref{eq:pplus}. For example, 
\begin{align*}
\alpha^{k}_k(x) =&\ 1 - \frac{11}{6} \Big(\frac{x-\delta k}{\delta} \Big) + \Big(\frac{x-\delta k}{\delta} \Big)^2 - \frac{1}{6}\Big(\frac{x-\delta k}{\delta} \Big)^3 - \frac{23}{3}\Big(\frac{x-\delta k}{\delta} \Big)^4 \notag \\
&+\frac{41}{2}\Big(\frac{x-\delta k}{\delta} \Big)^5 - \frac{55}{3}\Big(\frac{x-\delta k}{\delta} \Big)^6 + \frac{11}{2}\Big(\frac{x-\delta k}{\delta} \Big)^7.
\end{align*}
It is straightforward to check that 
\begin{align*}
\sum_{i=0}^{4} \alpha^{k}_{k+i}(x) = 1, \quad \alpha_{k}^{k}(\delta k) = 1, \quad \text{ and } \quad \alpha_{k+i}^{k} (\delta k) = 0.
\end{align*} 
The weights are degree-$7$ polynomials in $(x-\delta k)/ \delta$ whose coefficients  do not depend on $k$ or $\delta$, i.e., $\alpha^{k}_{k+i}(x)  = J_{i}\big((x-\delta k)/ \delta \big)$ 
for some polynomial $J_i(\cdot)$. Consequently, for any $x \in \R$, 
\begin{align*}
\alpha^{k+j}_{k+j+i}(x+ \delta j) = J_{i }\Big(\frac{x + \delta j- \delta(k+j)}{\delta} \Big)= J_{i}\Big(\frac{x-\delta k}{\delta} \Big) = \alpha^{k}_{k+i}(x),\quad j,k \in \Z,\ 0 \leq i \leq 4.
\end{align*} 
\finishproof

We now generalize Theorem~\ref{thm:1dinterpolant_def} and define an interpolation operator that can interpolate any function defined on   $ K \cap \delta \Z^{d}$ where $K \subset \R^{d}$ is convex.  The interpolator is based on forward differences, but one could also use central  or backward  differences to accommodate different domains shapes. The following theorem summarizes the key properties we want from it. 
\begin{theorem}
\label{thm:interpolant_def} 
Let the weights $\{\alpha_{k+i}^{k}: \R \to \R \ :\ k \in \Z,\ i = 0,1,2,3,4 \}$ be as in Theorem~\ref{thm:1dinterpolant_def} and suppose we are given a convex set $K \subset \R^{d}$ and a function  $f:  K \cap \delta \Z^{d} \to \R$. Letting $i = (i_1, \ldots, i_d) \in \Z^{d}$, we use the weights to define 
\begin{align}
A f(x) =&\  \sum_{i_d = 0}^{4} \alpha_{k_d(x)+i_d}^{k_d(x)}(x_d)\cdots \sum_{i_1 = 0}^{4} \alpha_{k_1(x)+i_1}^{k_1(x)}(x_1) f(\delta(k(x)+i)) \notag \\
=&\ \sum_{i_1, \ldots, i_d = 0}^{4} \bigg(\prod_{j=1}^{d} \alpha_{k_j(x) +i_j}^{k_j(x)   }(x_j)\bigg) f(\delta k(x)+i) , \quad x \in \text{Conv}(K_4),  \label{eq:af2}
\end{align}
where  $k(x) \in  \Z^{d}$ is defined by $k_{i}(x) = \lfloor x_i/\delta\rfloor$, and 
\begin{align*}
K_{4} = \{x \in K \cap \delta \Z^{d} : \delta(k(x)+ i) \in K \cap \delta \Z^{d} \text{ for all } 0 \leq i \leq 4e\}.
\end{align*} 
Then $A f(x) \in C^{3}(\text{Conv}(K_4))$   and is infinitely differentiable almost everywhere on $\text{Conv}(K_4) $. Also,
\begin{align}
A f(\delta k) = f(\delta k), \quad \delta k \in K \cap \delta \Z^{d}, \label{eq:interpolates2}
\end{align} 
and  there exists a constant $C > 0$ independent of $f(\cdot)$,$x$, and $\delta$, such that  
\begin{align}
\Big|  \frac{\partial^{a}}{\partial x^{a}} Af(x) \Big|   \leq&\  C \delta^{-\norm{a}_{1}}  \max_{\substack{ 0 \leq i_j \leq 4-a_j \\ j = 1,\ldots, d}} \abs{\Delta_{1}^{a_1}\ldots \Delta_{d}^{a_d} f(\delta (k(x)+i))}, \quad x \in \text{Conv}(K_4), \label{eq:multibound2}
\end{align}
for $0 \leq \norm{a}_{1} \leq 3$, and \eqref{eq:multibound2} also holds when $\norm{a}_{1} = 4$  for almost all $x \in  \text{Conv}(K_4) $.
\end{theorem}
Note that for any $J \subset \{1, \ldots, d\}$ and $J^{c} = \{1,\ldots, d\} \setminus J$, we may rewrite \eqref{eq:af2} as
\begin{align}
Af(x) =&\ \sum_{\substack{i_j = 0 \\ j \in J^{c}}}^{4} \Bigg( \bigg(\prod_{j\in J^{c}} \alpha_{k_j+i_j}^{k_j}(x_j)\bigg)  \sum_{\substack{i_j = 0 \\ j \in J}}^{4}\bigg(\prod_{j \in J} \alpha_{k_j+i_j}^{k_j}(x_j)\bigg) f(\delta(k+i))\Bigg). \label{eq:alpha2}
\end{align}
The representation in \eqref{eq:alpha2} will come in handy later on.  We define  
\begin{align}
F_k(x)=&\  \sum_{i_d = 0}^{4} \alpha_{k_d+i_d}^{k_d}(x_d)\cdots \sum_{i_1 = 0}^{4} \alpha_{k_1+i_1}^{k_1}(x_1) f(\delta(k+i)), \quad x \in \R^{d}, k \in K_4 \label{eq:alpha1}
\end{align}
to be the multidimensional analog of $P_{k}(x)$.
 Note that $A f(x)$ defined in Theorem~\ref{thm:interpolant_def} satisfies $A f(x) = F_{k(x)}(x)$ for $x \in \text{Conv}(K_4)$. Furthermore, \eqref{eq:alphas_interpolate} of Theorem~\ref{thm:1dinterpolant_def} implies \eqref{eq:interpolates2}. To prove Theorem~\ref{thm:interpolant_def}, it remains to verify the smoothness of $A f(x)$ and \eqref{eq:multibound2}.

 For any $x \in \R^{d}$ and $J \subset \{1, \ldots, d\}$, we write $x_{J}$ to denote the vector whose $i$th element equals $x_i 1(i \in J)$. The following result is the multidimensional analog of \eqref{eq:onedsmooth}. We prove it at the end of this section.
\begin{lemma}
\label{lem:multidinterp} 
Fix $k \in K_4$. For $u \in  [0,1]^{d}$, let $\Theta(u) = \{i: u_i = 1\}$ and $\Theta(u)^{c} = \{1,\ldots, d\}\setminus \Theta(u) $. Then 
\begin{align}
& \frac{\partial^{a}}{\partial x^{a}} F_k(x) \Big|_{x = \delta(k+u)} = \frac{\partial^{a}}{\partial x^{a}}F_{k+e_{\Theta(u)}}(x) \Big|_{x=\delta(k+u)} \quad \text{ for }  0 \leq \norm{a}_{1} \leq 3. \label{eq:multipaste}
\end{align}
Furthermore, there exists a constant $C> 0$ independent of $f(\cdot)$, $k$, and $\delta$  such that 
\begin{align}
 \Big| \frac{\partial^{a}}{\partial x^{a}}  F_k(x) \Big|  \leq&\  C \delta^{-\norm{a}_{1}}\bigg(\prod_{j=1}^{d}\Big(1 +  \Big|\frac{x_j-\delta k_j}{\delta} \Big| \Big)^{7-a_j} \bigg) \max_{\substack{ 0 \leq i_j \leq 4-a_j \\ j = 1,\ldots, d}} \abs{\Delta_{1}^{a_1}\ldots \Delta_{d}^{a_d} f(\delta (k+i))} \label{eq:premultibound}
\end{align}
 for all $0 \leq \norm{a}_{1} \leq 4$ and all $x \in \text{Conv}(K_4)$ where the derivative above is well defined.
 
\end{lemma}
\noindent 
Lemma~\ref{lem:multidinterp} implies Theorem~\ref{thm:interpolant_def}. Indeed, \eqref{eq:multipaste} implies $A f(x) \in C^{3}(\text{Conv}(K_4))$, and since $\alpha_{k+i}^{k}(x) \in C^{\infty}(\R)$, we know $A f(x)$ is infinitely differentiable everywhere except at the points where the  $F_k(x)$ are glued together, i.e.,  on the set $\{x \in \text{Conv}(K_4)\ |\  x_i \in \delta \Z \text{ for some } i \in \{1, \ldots, d\} \}$, which has Lebesgue measure zero.   Furthermore, \eqref{eq:multibound2} follows directly from \eqref{eq:premultibound}. We now prove Lemma~\ref{lem:multidinterp}.
 
\startproof{Proof of Lemma~\ref{lem:multidinterp}}
We first prove \eqref{eq:multipaste}. Fix $k \in K_{4}$ and let $j' \in \Theta(u)$. From \eqref{eq:alpha2} it follows that 
\begin{align*}
 \frac{\partial^{a}}{\partial x^{a}} F_k(x)\Big|_{x = \delta(k+u)} =&\  \sum_{\substack{i_j = 0 \\ j \in \{1, \ldots, d\}\setminus \{j'\}}}^{4} \Bigg( \bigg(\prod_{j \in \{1, \ldots, d\}\setminus \{j'\}} \frac{\partial^{a_j}}{\partial x_{j}^{a_j}} \alpha_{k_j+i_j}^{k_j}(x_j)\Big|_{x_j = \delta(k_j+u_j)}\bigg) \\
& \hspace{2cm} \times  \sum_{ i_{j'} = 0 }^{4}\bigg(  \frac{\partial^{a_{j'}}}{\partial x_{j'}^{a_{j'}}}\alpha_{k_{j'}+i_{j'}}^{k_{j'}}(x_{j'})\Big|_{x_{j'} = \delta(k_{j'}+1)} \bigg) f(\delta(k+i))\Bigg).
\end{align*}
For the inner sum, note that 
\begin{align*}
& \sum_{ i_{j'} = 0 }^{4}\bigg(  \frac{\partial^{a_{j'}}}{\partial x_{j'}^{a_{j'}}}\alpha_{k_{j'}+i_{j'}}^{k_{j'}}(x_{j'})\Big|_{x_{j'} = \delta(k_{j'}+1)} \bigg) f(\delta(k+i)) \\
 =&\  \sum_{ i_{j'} = 0 }^{4}\bigg(\frac{\partial^{a_{j'}}}{\partial x_{j'}^{a_{j'}}}\alpha_{k_{j'}+1+i_{j'}}^{k_{j'}+1}(x_{j'})\Big|_{x_{j'} = \delta(k_{j'}+1)} \bigg) f(\delta(k+i + e_{j'})),
\end{align*}
which follows from \eqref{eq:onedsmooth}. Repeating the above procedure for all other elements of $\Theta(u)$, we see that 
\begin{align*}
 \frac{\partial^{a}}{\partial x^{a}} F_k(x)\Big|_{x = \delta(k+u)} =&\ \sum_{\substack{i_j = 0 \\ j \in \Theta(u)^{c}}}^{4} \Bigg( \bigg(\prod_{j\in \Theta(u)^{c}} \frac{\partial^{a_j}}{\partial x_{j}^{a_j}} \alpha_{k_j+i_j}^{k_j}(x_j)\Big|_{x_j = \delta(k_j+u_j)}\bigg) \\
& \hspace{1cm} \times   \sum_{\substack{i_j = 0 \\ j \in \Theta(u)}}^{4}\bigg(\prod_{j \in \Theta(u)} \frac{\partial^{a_j}}{\partial x_{j}^{a_j}}\alpha_{k_j+1+i_j}^{k_j+1}(x_j)\Big|_{x_j = \delta(k_j+1)} \bigg) f(\delta(k+i+e_{\Theta(u)}))\Bigg)\\
=&\ \frac{\partial^{a}}{\partial x^{a}} F_{k+e_{\Theta(u)}}(x) \Big|_{x=\delta(k+u)},
\end{align*}  
which proves \eqref{eq:multipaste}. It remains to prove the bound on $\big| \frac{\partial^{a}}{\partial x^{a}} F_k(x)\big| $ in \eqref{eq:premultibound}. We know
\begin{align*}
\frac{\partial^{a}}{\partial x^{a}}F_k(x) =&\ \sum_{i_1=0}^{4}\frac{\partial^{a_1}}{\partial x_{1}^{a_1}}  \alpha_{k_1+i_1}^{k_1}(x_1)   \cdots\sum_{i_d=0}^{4}\frac{\partial^{a_d}}{\partial x_{d}^{a_d}}  \alpha_{k_d+i_d}^{k_d}(x_d)  f(\delta(k+i)).
\end{align*}
%Let us treat $i_1, \ldots, i_{d-1}$ as fixed and consider $f(\delta (k+i))$ a one-dimensional function of $\delta (k_d+i_d)$, i.e. 
%\begin{align*}
% \sum_{i_2=0}^{4}\frac{\partial^{a_2}}{\partial x_{2}^{a_2}}  \alpha_{k_2+i_2}^{k_2}(x_2)  \cdots\sum_{i_d=0}^{4}\frac{\partial^{a_d}}{\partial x_{d}^{a_d}}  \alpha_{k_d+i_d}^{k_d}(x_d)  f(\delta(k+i)) = \bar f(\delta (k_1+i_1)).
%\end{align*}
 By inspecting the form of the one-dimensional  $P_k(\cdot)$ in \eqref{eq:pplus}, one can check that 
\begin{align*}
\sum_{i_d=0}^{4}\frac{\partial^{a_d}}{\partial x_{d}^{a_d}}  \alpha_{k_d+i_d}^{k_d}(x_d)    f(\delta (k +i )) = \delta^{-a_{d}} Q^{(d)}\Big( \frac{x_d - \delta k_d}{\delta} \Big) ,
\end{align*}
where $Q^{(d)}(\cdot)$ is a degree-$(7-a_{d})$  polynomial whose coefficients depend on $ f(\cdot)$ only through 
\begin{align*}
&\Delta_{d}^{a_d} f\big(\delta (k + (i_1,\ldots, i_{d-1},i_d))\big),  \quad \text{ for }  i_d = 0, \ldots, 4-a_{d},
\end{align*}
and are independent of $\delta$. This implies in particular that 
\begin{align*}
 \Big| \sum_{i_d=0}^{4}\frac{\partial^{a_d}}{\partial x_{d}^{a_d}} \alpha_{k_d+i_d}^{k_d}(x_d)  f(\delta(k+i)) \Big|  \leq&\ C \delta^{-a_{d}}  \Big(1 +  \Big|\frac{x_d-\delta k_d}{\delta} \Big| \Big)^{7-a_{d}}  \max_{\substack{ 0 \leq i_d \leq 4-a_{d} }} \abs{\Delta_{d}^{a_d} f(\delta (k+i))}.  
\end{align*}
We now consider 
\begin{align*}
\sum_{i_{d-1}=0}^{4}\frac{\partial^{a_{d-1}}}{\partial x_{{d-1}}^{a_{d-1}}}  \alpha_{k_{d-1}+i_{d-1}}^{k_{d-1}}(x_{d-1}) \bigg( \sum_{i_d=0}^{4}\frac{\partial^{a_d}}{\partial x_{d}^{a_d}}  \alpha_{k_d+i_d}^{k_d}(x_d)  f(\delta(k+i)) \bigg).
\end{align*}
When viewed as a one-dimensional function of $x_{d-1}$, the above is again a degree-$(7-a_{d-1})$ polynomial that depends on the quantity inside the parentheses only through 
\begin{align*}
&\Delta_{d-1}^{a_{d-1}}\bigg( \sum_{i_d=0}^{4}\frac{\partial^{a_d}}{\partial x_{d}^{a_d}}  \alpha_{k_d+i_d}^{k_d}(x_d)  f(\delta(k+i)) \bigg), \quad \text{ for } i_{d-1} = 0, \ldots, 4-a_{d-1}.
\end{align*}
Hence, 
\begin{align*}
&\bigg| \sum_{i_{d-1}=0}^{4}\frac{\partial^{a_{d-1}}}{\partial x_{{d-1}}^{a_{d-1}}}  \alpha_{k_{d-1}+i_{d-1}}^{k_{d-1}}(x_{d-1}) \bigg( \sum_{i_d=0}^{4}\frac{\partial^{a_d}}{\partial x_{d}^{a_d}}  \alpha_{k_d+i_d}^{k_d}(x_d)  f(\delta(k+i)) \bigg) \bigg| \\
 \leq&\ C \delta^{-a_{d-1}-a_{d}} \Big(1 +  \Big|\frac{x_{d-1}-\delta k_{d-1}}{\delta} \Big| \Big)^{7-a_{d-1}}  \Big(1 +  \Big|\frac{x_d-\delta k_d}{\delta} \Big| \Big)^{7-a_{d}}  \\
 & \hspace{4.5cm} \times \max_{\substack{ 0 \leq i_d \leq 4-a_{d} \\ 0 \leq i_{d-1} \leq 4 - a_{d-1} }} \abs{\Delta_{d-1}^{a_{d-1}} \Delta_{d}^{a_d} f(\delta (k+i))}.
\end{align*}
Repeating this argument along each of the remaining $d-2$ dimensions proves \eqref{eq:premultibound}. 
\finishproof

\section{Interchange in Multiple Dimensions} 
\label{sec:mdimunbound}
In this section we prove Proposition~\ref{lem:interchanged1} by proving the more general Proposition~\ref{lem:mdiminterchange} stated below.  Consider a CTMC living on $\Z^{d}$ with generator 
\begin{align*}
G_X f(\delta k) =&\ \sum_{\ell \in \Z^{d}} \beta_{\ell}(\delta k) (f(\delta (k+\ell))-f(\delta k)), \quad k \in \Z^{d}.
\end{align*} 

\begin{proposition} \label{lem:mdiminterchange}
Fix $f: \delta \Z^{d} \to \R$ and assume  that
\begin{align}
 \sum_{\ell \in \Z^{d}} \abs{\beta_{\ell}(\delta k) (f(\delta (k+\ell))-f(\delta k))} < \infty, \quad k \in \Z^{d}, \label{eq:intergrab2}
\end{align}
which is trivially satisfied  when the number of possible transitions from each state is finite. For $x \in \R^{d}$ define  $k(x) \in  \Z^{d}$   by $k_{i}(x) = \lfloor x_i/\delta\rfloor $. Then
\begin{align}
 A G_{X} f(x) =&\  \sum_{\ell \in \Z}  A \beta_{\ell}(x)  \big( A f(x+\delta \ell) - A f(x)\big) + \varepsilon(x), \quad x \in \R^{d},    \label{eq:intererror2}
\end{align}
 where
\begin{align}
\varepsilon(x) =&\  \sum_{\ell \in \Z} \sum_{i_1, \ldots, i_d = 0}^{4} \bigg(\prod_{j=1}^{d} \alpha_{k_j(x) +i_j}^{k_j(x)   }(x_j)\bigg) \Big(\beta_{\ell}\big(\delta (k(x)+i)\big) - A \beta_{\ell}(x)\Big) \notag \\
& \hspace{3cm} \times \Big(f\big(\delta (k(x)+\ell+i)\big)-f\big(\delta (k(x)+i)\big) - \big(f\big(\delta (k(x)+\ell )\big)-f (\delta k(x) )\big)\Big). \label{eq:varepsmdim}
\end{align}
%\begin{align*}
%\varepsilon(x) =&\ \sum_{\ell \in \Z}  \sum_{i=0}^{4} \alpha^{k(x)}_{k(x)+i}(x)\Big(\beta_{\ell}\big(\delta (k(x)+i)\big) - A \beta_{\ell}(x)\Big) \notag  \\
%& \hspace{1.5cm} \times \Big(f\big(\delta (k(x)+\ell+i)\big)-f\big(\delta (k(x)+i)\big) - \big(f\big(\delta (k(x)+\ell )\big)-f\big(\delta k(x) \big)\big)\Big).
%\end{align*}
\end{proposition}
Before proving Proposition~\ref{lem:mdiminterchange}, let us reconcile the forms of  $\varepsilon(x)$ in \eqref{eq:varepsmdim} above  and in \eqref{eq:vareps} of Proposition~\ref{lem:interchanged1}. When $d = 1$, \eqref{eq:varepsmdim} equals 
\begin{align*}
 &\sum_{\ell \in \Z}  \sum_{i=0}^{4} \alpha^{k(x)}_{k(x)+i}(x)\Big(\beta_{\ell}\big(\delta (k(x)+i)\big) - A \beta_{\ell}(x)\Big)  \\
& \hspace{3cm} \times \Big(f\big(\delta (k(x)+\ell+i)\big)-f\big(\delta (k(x)+i)\big) - \big(f\big(\delta (k(x)+\ell )\big)-f\big(\delta k(x)\big)\big)\Big).
\end{align*}
Using a telescoping series, we see that if $\ell > 0$, 
\begin{align*}
& f\big(\delta (k(x)+\ell+i)\big)-f\big(\delta (k(x)+i)\big) - \big(f\big(\delta (k(x)+\ell )\big)-f\big(\delta k(x)\big)\big)\\
=&\ \sum_{j=0}^{\ell-1} \big(\Delta f(\delta(k(x)+j+\ell)) -  \Delta f(\delta(k(x) +j))\big) \\
=&\ \sum_{j=0}^{i-1} \sum_{m=0}^{\ell-1}   \Delta^2 f(\delta(k(x)+m +j)).
\end{align*} 
Similarly, when $\ell < 0$,  
\begin{align*}
& f\big(\delta (k(x)+\ell+i)\big)-f\big(\delta (k(x)+i)\big) - \big(f\big(\delta (k(x)+\ell )\big)-f\big(\delta k(x)\big)\big)\\
  =&\ - \sum_{j=0}^{i-1} \sum_{m=-\ell}^{ -1}  \Delta^2 f(\delta(k(x)+m +j)).
\end{align*} 
Therefore,  Propositions~\ref{lem:interchanged1} and \ref{lem:mdiminterchange} are equivalent when $d = 1$.  When $d > 1$, it is also possible to write \eqref{eq:varepsmdim} as a telescoping series of second-order differences of $f(\delta k)$. We leave this as an exercise to the interested reader.

\startproof{Proof of Proposition~\ref{lem:mdiminterchange}}
Fix $x \in \R^{d}$. We will write $k$ instead of $k(x)$ for convenience. Recalling the form of $A$ from Theorem~\ref{thm:interpolant_def},
%that for any function $f: \delta \Z^{d} \to \R$, 
%\begin{align*}
%Af(x) = \sum_{i_1, \ldots, i_d = 0}^{4} \bigg(\prod_{j=1}^{d} \alpha_{k_j +i_j}^{k_j   }(x_j)\bigg) f\big(\delta(k+i)\big).
%\end{align*} 
it follows that $A G_{X} f(x)$ equals
\begin{align}
&  \sum_{i_1, \ldots, i_d = 0}^{4} \bigg(\prod_{j=1}^{d} \alpha_{k_j +i_j}^{k_j   }(x_j)\bigg) \sum_{\ell \in \Z} \beta_{\ell}\big(\delta (k+i)\big) \Big(f\big(\delta (k+\ell+i)\big)-f\big(\delta (k+i)\big)\Big) \notag \\
=&\  \sum_{\ell \in \Z}  A \beta_{\ell}(x) \sum_{i_1, \ldots, i_d = 0}^{4} \bigg(\prod_{j=1}^{d} \alpha_{k_j +i_j}^{k_j   }(x_j)\bigg)  \Big(f\big(\delta (k+\ell+i)\big)-f\big(\delta (k+i)\big)\Big)  \label{eq:first} \\
&+ \sum_{\ell \in \Z} \sum_{i_1, \ldots, i_d = 0}^{4} \bigg(\prod_{j=1}^{d} \alpha_{k_j +i_j}^{k_j   }(x_j)\bigg)\Big(\beta_{\ell}\big(\delta (k+i)\big) - A \beta_{\ell}(x)\Big)   \Big(f\big(\delta (k+\ell+i)\big)-f\big(\delta (k+i)\big)\Big). \label{eq:second}
\end{align}
Interchanging the summations is allowed by the Fubini-Tonelli theorem due to assumption \eqref{eq:intergrab2}. Consider first the inner sum in \eqref{eq:first} and observe that for each $\ell \in \Z^{d}$, 
\begin{align*}
& \sum_{i_1, \ldots, i_d = 0}^{4} \bigg(\prod_{j=1}^{d} \alpha_{k_j +i_j}^{k_j   }(x_j)\bigg)  f\big(\delta (k+\ell+i)\big)- \sum_{i_1, \ldots, i_d = 0}^{4} \bigg(\prod_{j=1}^{d} \alpha_{k_j +i_j}^{k_j   }(x_j)\bigg) f\big(\delta (k+i)\big) \\
=&\ \sum_{i_1,\ldots,i_d = 0}^{4}\bigg(\prod_{j=1}^{d} \alpha_{k_j+\ell_{j}  +i_j}^{k_j+\ell_{j}}(x_j+\delta \ell_{j})\bigg)   f\big(\delta (k+\ell+i)\big)- A f(x)\\
=&\ A f(x+\delta \ell) - A f(x), 
\end{align*} 
where in the first equality we used the translation invariance property of the weights stated in  \eqref{eq:weights} of Theorem~\ref{thm:1dinterpolant_def}. Moving on, we see that \eqref{eq:second} equals
\begin{align*}
 &  \sum_{\ell \in \Z} \sum_{i_1, \ldots, i_d = 0}^{4} \bigg(\prod_{j=1}^{d} \alpha_{k_j +i_j}^{k_j   }(x_j)\bigg) \Big(\beta_{\ell}\big(\delta (k+i)\big) - A \beta_{\ell}(x)\Big)  \\
& \hspace{3cm} \times \Big(f\big(\delta (k+\ell+i)\big)-f\big(\delta (k+i)\big) - \big(f\big(\delta (k+\ell )\big)-f (\delta k )\big)\Big)\\
&+ \sum_{\ell \in \Z}  \Big(f\big(\delta (k+\ell )\big)-f (\delta k )\Big)\sum_{i_1, \ldots, i_d = 0}^{4} \bigg(\prod_{j=1}^{d} \alpha_{k_j +i_j}^{k_j   }(x_j)\bigg) \Big(\beta_{\ell}\big(\delta (k+i)\big) - A \beta_{\ell}(x)\Big)  .
\end{align*}
%\begin{align*}
%&A G_{X} f(x) \\
%=&\ \sum_{\ell \in \Z}  \sum_{i=0}^{4} \alpha^{k}_{k+i}(x) \beta_{\ell}\big(\delta (k+i)\big) \Big(f\big(\delta (k+\ell+i)\big)-f\big(\delta (k+i)\big)\Big)\\
%=&\  \sum_{\ell \in \Z}  A \beta_{\ell}(x) \sum_{i=0}^{4} \alpha^{k}_{k+i}(x)   \Big(f\big(\delta (k+\ell+i)\big)-f\big(\delta (k+i)\big)\Big) \\
%&+ \sum_{\ell \in \Z}  \sum_{i=0}^{4} \alpha^{k}_{k+i}(x)\Big(\beta_{\ell}\big(\delta (k+i)\big) - A \beta_{\ell}(x)\Big)   \Big(f\big(\delta (k+\ell )\big)-f\big(\delta (k )\big)\Big)\\
%&+ \sum_{\ell \in \Z}  \sum_{i=0}^{4} \alpha^{k}_{k+i}(x)\Big(\beta_{\ell}\big(\delta (k+i)\big) - A \beta_{\ell}(x)\Big)  \\
%& \hspace{3cm} \times \Big(f\big(\delta (k+\ell+i)\big)-f\big(\delta (k+i)\big) - \big(f\big(\delta (k+\ell )\big)-f\big(\delta (k )\big)\big)\Big).
%\end{align*}
The second line  equals zero because \eqref{eq:weights_sum_one} of Theorem~\ref{thm:1dinterpolant_def}  implies $\sum_{i_1, \ldots, i_d = 0}^{4} \Big(\prod_{j=1}^{d} \alpha_{k_j +i_j}^{k_j   }(x_j)\Big) = 1$   and because $A \beta_{\ell}(x) = \sum_{i_1, \ldots, i_d = 0}^{4} \Big(\prod_{j=1}^{d} \alpha_{k_j +i_j}^{k_j   }(x_j)\Big) \beta_{\ell}\big(\delta (k+i)\big)$  by definition.  
\finishproof

\section{Proofs of Miscellaneous Technical Lemmas}

\subsection{Proof of Lemma~\ref{lem:convdet}}
\label{sec:proofconvdet}

We now state and prove an auxiliary result, and then prove Lemma~\ref{lem:convdet}.

\begin{lemma} \label{lem:auxconv}
Suppose $g^*:\R \to \R$ is three times continuously differentiable with an absolutely continuous third derivative,  and let $g: \delta \Z \to \R$ be its restriction to $\delta \Z $.  For $1 \leq v \leq 4$, 
\begin{align}
\Delta^{v}  g(\delta k) = \delta^{v} \int_{0}^{v \delta} c^{v}(u)  \frac{\partial^{v}}{\partial x^{v}} g^*(\delta k + u) du,   \label{eq:cver}
\end{align}
where $c^{v}(x)$ is a function such that $\sup_{x \in [0,\delta v]} \abs{c^{v}(x)} \leq C $ for some constant $C > 0$ independent of $k$, $g^*(x)$, and $\delta$.   
\end{lemma}

\startproof{Proof of Lemma~\ref{lem:auxconv}}
Suppose  $v = 4$; the case $1 \leq v \leq 3$ is handled similarly. Using Taylor expansion, we have 
 \begin{align*}
\Delta^{4}  g(\delta k) =&\  \Delta^{3} \big(g (\delta (k+1)) - g(\delta k) \big)\\ 
=&\ \Delta^{3} \Big( \delta (g^*)'(\delta k) + \frac{1}{2}\delta^{2} (g^*)''(\delta k) + \frac{1}{6} \delta^{3} (g^*)'''(\delta k) + \frac{1}{6}\int_{0}^{\delta} g^{(4)}(\delta k + u) (\delta - u)^{3} du  \Big).
\end{align*}
We note that
\begin{align*}
&\Delta^{3} \Big( \frac{1}{6}\int_{0}^{\delta} g^{(4)}(\delta k + u) (\delta - u)^{3} du  \Big) \\
=&\ \frac{1}{6}\Bigg( \int_{0}^{\delta} (g^*)^{(4)}(\delta (k+3) + u) (\delta - u)^{3} du - 3 \int_{0}^{\delta} (g^*)^{(4)}(\delta (k+1) + u) (\delta - u)^{3} du \\
& \qquad + 3\int_{0}^{\delta} (g^*)^{(4)}(\delta (k+2) + u) (\delta - u)^{3} du - \int_{0}^{\delta} (g^*)^{(4)}(\delta k + u) (\delta - u)^{3} du \Bigg).
\end{align*}
Applying a similar expansion to $\Delta^{3} \big( \delta (g^*)'(\delta k) + \frac{1}{2}\delta^{2} (g^*)''(\delta k) + \frac{1}{6} \delta^{3} (g^*)'''(\delta k)  \big)$ proves \eqref{eq:cver}.
\finishproof

\startproof{Proof of Lemma~\ref{lem:convdet}}
Since $h^*(X) = h(X)$, the triangle inequality implies that
\begin{align*}
\big| \E h^*(X) - \E h^*(Y) \big|  \leq&\ \big| \E h(X ) - \E Ah( Y) \big| + \big| \E Ah( Y) - \E h^*( Y) \big|.
\end{align*}
For $x \in \R^{d}$ let $k(x) \in  \Z^{d}$ be defined  by $k_{i}(x) = \lfloor x_i/\delta\rfloor$. Since $h^*(\delta k(x)) = h(\delta k(x)) = A h(\delta k(x))$ and $\abs{x_i - k_i(x)} \leq \delta$,  
\begin{align*}
\big| \E Ah(Y) - \E h^*(Y) \big| \leq&\ \big| \E Ah(Y) - \E h^*(\delta k(Y)) \big| +  \big| \E h^*(\delta k(Y)) - \E h^*(Y) \big|\\
=&\ \big| \E Ah(Y) - \E Ah(\delta k(Y)) \big| +  \big| \E h^*(\delta k(Y)) - \E h^*(Y) \big|\\
\leq&\ C \delta \max_{\substack{  1 \leq j \leq d  }}\sup_{\substack{   x \in \R^{d} }} \bigg| \frac{\partial}{\partial x_{j}} Ah(x)  \bigg| + C \delta \max_{\substack{  1 \leq j \leq d  }}\sup_{\substack{   x \in \R^{d} }} \bigg| \frac{\partial}{\partial x_{j}} h^{*}(x)  \bigg|.
\end{align*}
  Using the bound in  \eqref{eq:multibound2} from Theorem~\ref{thm:interpolant_def}, it follows that 
\begin{align*}
\max_{\substack{  1 \leq j \leq d  }}\sup_{\substack{   x \in \R^{d} }} \bigg| \frac{\partial}{\partial x_{j}} Ah(x)\bigg|  \leq   C \delta^{-1} \max_{\substack{  1 \leq j \leq d  }}\sup_{\substack{   x \in \R^{d} }} \abs{ \Delta_{j}  h(\delta k(x))} =&\   C \delta^{-1}\max_{\substack{  1 \leq j \leq d  }}\sup_{\substack{   x \in \R^{d} }} \abs{ \Delta_{j}  h^{*}(\delta k(x))} \\
\leq&\   C   \max_{\substack{  1 \leq j \leq d  }}\sup_{\substack{   x \in \R^{d} }} \bigg| \frac{\partial}{\partial x_{j}} h^{*}(x)  \bigg|.
\end{align*}  
This proves the first claim. The other two claims follow by observing that if $h^{*} \in \lipone$, then the mean-value theorem implies $h \in \dlipone$, and if $h^{*} \in \mathcal{M}$, then we can apply \eqref{eq:cver} along each dimension to show that $h  \in \mathcal{M}_{disc}(C')$ for some $C' > 0$. 
\finishproof

\subsection{Proof of Lemma~\ref{lem:genpoisson}}
\label{sec:proofgenpoisson}
\startproof{Proof of Lemma~\ref{lem:genpoisson}} 
Let $\tau = \inf_{t \geq 0} \{X(t) \neq X(0)\}$ be the first jump time. For any $\varepsilon > 0$, 
\begin{align*}
g(\delta k) =&\ \int_{0}^{\infty} \E_{X(0) = \delta k}  \Big(\big(h(X(t)) - \E h(X) \big) 1(t <  \tau \wedge \varepsilon)\Big)   dt  \\
&+ \int_{0}^{\infty} \E_{X(0) = \delta k}  \Big(\big(h(X(t)) - \E h(X) \big) 1(t > \tau \wedge \varepsilon)\Big)   dt.
\end{align*}
Set $r = \sum_{k' \neq k } q_{k,k'}$ and note that $\Prob(\tau > \varepsilon) = e^{-r \varepsilon}$. The strong Markov property implies
\begin{align*}
\int_{0}^{\infty} \E_{X(0) = \delta k}  \Big(\big(h(X(t)) - \E h(X) \big) 1(t > \tau \wedge \varepsilon)\Big)   dt =&\ e^{-r \varepsilon} g(\delta k) + (1- e^{-r \varepsilon}) \sum_{k' \neq k } \frac{q_{k,k'}}{r} g(\delta k').
\end{align*}
Furthermore, since $X(t) = X(0) $ for $t < \tau \wedge \varepsilon$, 
\begin{align*}
\int_{0}^{\infty} \E_{X(0) = \delta k}  \Big(\big(h(X(t)) - \E h(X) \big) 1(t <  \tau \wedge \varepsilon)\Big)   dt =&\ \big(h(\delta k) - \E h(X)\big) \int_{0}^{\infty} \Prob_{X(0) = \delta k} (\tau \wedge \varepsilon > t) dt\\
=&\ \big(h(\delta k) - \E h(X)\big) \E \big(\tau \wedge \varepsilon \big).
\end{align*}
Combining the three equations above and rearranging terms yields
\begin{align*}
0=&\  \big(h(\delta k) - \E h(X)\big) \E \big(\tau \wedge \varepsilon \big)   + (1- e^{-r \varepsilon}) \sum_{k' \neq k } \frac{q_{k,k'}}{r} \big( g(\delta k') - g(\delta k)\big).
\end{align*}
We conclude by dividing both sides by $\varepsilon$ and taking $\varepsilon \to 0$, and by noting that 
\begin{align*}
\frac{1}{\varepsilon} \E \big(\tau \wedge \varepsilon \big) = e^{-r \varepsilon} + (1 - e^{-r \varepsilon}) \frac{1}{\varepsilon}\E \big((\tau \wedge \varepsilon) 1(\tau < \varepsilon) \big) \to 1, \quad \text{ as } \varepsilon \to 0.
\end{align*}
\finishproof

 \end{APPENDICES}

% Acknowledgments here
\ACKNOWLEDGMENT{The author would like to thank Han Liang Gan for stimulating discussions during early stages of this work, as well as Robert Bray and Shane Henderson for providing feedback on early drafts. The author is also grateful to two anonymous referees for their numerous suggestions to improve the presentation of the material, as well as Zhe Su whose input helped significantly reduce the length of the manuscript. }

% References here (outcomment the appropriate case)

% CASE 1: BiBTeX used to constantly update the references
%   (while the paper is being written).
\bibliographystyle{informs2014} % outcomment this and next line in Case 1
\bibliography{dai20190911} % if more than one, comma separated

%%%%%%%%%%%%%%%%%
\end{document}